\documentclass[10pt,a4paper]{article}
\usepackage[utf8]{inputenc}
\usepackage{natbib}
\usepackage[english]{babel}
\usepackage{fullpage}

\usepackage{amsfonts,amssymb,amsmath,amscd,latexsym,dsfont,amsthm}
\usepackage{graphicx}
\usepackage{natbib}
\usepackage{color}
\usepackage{enumitem}
\usepackage{url}

\newcounter{th} 

\newtheorem{ass}{Assumption}
\newtheorem{theorem}[th]{Theorem}
\usepackage{xr}
\usepackage{zref-base}
\usepackage{atveryend}
\makeatletter
\AfterLastShipout{%
  \if@filesw   
    \begingroup
      \zref@setcurrent{default}{\the\value{equation}}%
      \zref@wrapper@immediate{%
        \zref@labelbyprops{eq-last}{default}%
      }%
    \endgroup
  \fi
}
\makeatother

\usepackage{authblk}
\title{Full Model Estimation for Non-Parametric Multivariate Finite Mixture Models}
\author[1]{Marie Du Roy de Chaumaray}
\author[1]{Matthieu Marbac}
\affil[1]{Univ.  Rennes, Ensai, CNRS, CREST - UMR 9194, F-35000 Rennes, France}

\DeclareMathOperator*{\argmax}{arg\,max}

\newcommand{\simplex}{\mathcal{S}}
\newcommand{\bX}{\boldsymbol{X}}
\newcommand{\bx}{\boldsymbol{x}}
\newcommand{\bpi}{\boldsymbol{\pi}}
\newcommand{\bpsi}{\boldsymbol{\psi}}

\begin{document}
\maketitle

\begin{abstract}
This paper addresses the problem of full model estimation for non-parametric finite mixture models. It presents an approach for selecting the number of components and the subset of discriminative variables (\emph{i.e.,} the subset of variables having  different distributions among the mixture components). 
The proposed approach considers a discretization of each variable into $B$ bins  and a penalization of the resulting log-likelihood. 
Considering that the number of bins tends to infinity as the sample size tends to infinity, we prove that our estimator of the model (number of components and subset of relevant variables for clustering) is consistent under a suitable choice of the penalty term. Interest of our proposal is illustrated on simulated and benchmark data.
\end{abstract}


\section{Introduction}
Finite mixture models permits to achieve clustering by estimating the distribution of the observed variables \citep{McL00,Fruhwirth2019handbook}. This paper focuses on a full model selection (\emph{i.e.,} estimation of the number of components and detection of the subset of the relevant variables for clustering) for non-parametric mixture models where no assumptions are made on the component distribution except to be defined as a product of univariate densities (see \citet{ChauveauSurveys2015} for a review). 
Thus, we consider a sample composed of $n$ independent observations $\bX_1,\ldots,\bX_n$ where $\bX_i=(X_{i1},\ldots,X_{iJ})^\top\in\mathcal{X}$ is the vector composed of the $J$ variables collected on subject $i$ defined over the space $\mathcal{X}=\mathcal{X}_1\times \ldots\times\mathcal{X}_J$ where each $\mathcal{X}_j$ is compact. Each $\bX_i$ is identically distributed according to the non-parametric mixture of $K$ components defined by the density
\begin{equation}\label{eq:mixturegrl}
g(\bx_i)=\sum_{k=1}^K \pi_k \prod_{j=1}^J \eta_{kj}(x_{ij}) ,
\end{equation}
where  $\bpi=(\pi_1,\dots,\pi_K)^\top$ is a finite dimensional parameter belonging to the simplex of size $K$, $\simplex_K= \{\boldsymbol{u}\in[0,1]^K : \;\sum_{k=1}^K u_k=1\}$ and where  the univariate densities $\eta_{kj}$ constitute infinite dimensional parameters and are supposed to be striclty positive except on a set of lebesgue measure zero. Model \eqref{eq:mixturegrl} has been used in different fields like in behavioral science \citep{CloggBook1995}, econometry \citep{HuJofEco2013,CompianiEcoJ2016} or sociology \citep{HagenaarsBook2002}. One standard situation where the conditional independence assumption implied by \eqref{eq:mixturegrl} holds true, is the framework of standard repeated measure random effect model, where the suject-level effect is replaced by a component-level effect. 
 Among the recent developments related to \eqref{eq:mixturegrl}, one can cite the papers of   \citet{HallAOS2003}, \citet{HallBiometrika2005} and \citet{AllmanAOS2009} who studied the model identifiability, while \citet{BenagliaJCGS2009}, \citet{LevineBiometrika2011} and \citet{ZhengJASA2020} proposed  different algorithms for estimating the parameters  when the number of components $K$ is known.

Even if  the estimation of the groups is a key-point for clustering, estimating the number of components in \eqref{eq:mixturegrl} was still an open problem until recently. 
Thus, before the works of \citet{KasaharaJRSSB2014}, \citet{BonhommeJRSSB2016}, \citet{BonhommeAOS2016} and \citet{Kwon2019estimation}, no tools were available for selecting the number of components $K$. This was the main drawback to  non-parametric mixture models compared to parametric mixture models which allow the use of information criteria for model selection.
The two recent papers of \citet{KasaharaJRSSB2014} and \citet{Kwon2019estimation}  introduced two approaches for determining the smallest value of $K$ such that \eqref{eq:mixturegrl} holds true. Both works are based on the results of \citet{AllmanAOS2009} that established the model identifiability, using a theorem of \citet{KruskalLA1977}, whenever there are at least three variables with linearly independent univariate density functions  (except possibly on a set of Lebesgue measure zero).
\citet{KasaharaJRSSB2014} provide an estimation of the lower bound of the number of components by considering a partition of the support of each variable (\emph{e.g.,} using a decomposition into bins) and by using the identifiability of the latent class model (\emph{i.e.,} mixture models where each component is a product of multinomial distributions). This discretization allows to consider the tensor defining the probability of each event while the rank of this tensor permits to derive a lower bound on $K$. Note that previous works on non-parametric mixture models considered a bin decomposition (\emph{i.e.,} a specific discretization method) to estimate non-parametric mixture models or study their identifiability but not for model selection  \citep{HettmanspergerJRSSB2000,CruzJRSSC2004,ElmoreComStat2004}. However, \citet{KasaharaJRSSB2014} do not provide a method for selecting the discretization (\emph{i.e.,} number of elements, locations of those elements). Thus, their method is only consistent to a lower bound of $K$ (see Section 2.3 in \citet{Kwon2019estimation}). Alternatively, \citet{Kwon2019estimation} consider an integral operator, identified from the distribution of $X$, that has a rank equal to $K$. Noting that the singular values of operators are stable under perturbations (to handle the fact that this operator is estimated from the observed sample), a thresholding rule allowing to count the number of non-zero singular values provides a consistent estimator of $K$. One advantage of the approach of \citet{Kwon2019estimation} is to avoid the use of discretization, even if some connexions can be established with the approach of \citet{KasaharaJRSSB2014} (see Section 2.3 in \citet{Kwon2019estimation}). 
One elegant property of the methods of \citet{KasaharaJRSSB2014} and  \citet{Kwon2019estimation} is that both methods determine an estimator of $K$ without performing the density estimation for different numbers of clusters and without determining ahead a maximum number of clusters. Thus, those methods start with a step of model selection followed by the estimation of the parameters for the selected model. This is quite unusual because, when model selection is conducted for parametric mixture models via an information criterion, parameter estimation needs to be first performed for each competing model in order to compute the information criterion. Note that the use of the identifiability results \citep{AllmanAOS2009} is crucial to study the rank of the objects considered by \citet{KasaharaJRSSB2014} and \citet{Kwon2019estimation}.
The approaches of \citet{KasaharaJRSSB2014} and \citet{Kwon2019estimation} are mainly based on the distribution of a couple of variables. Thus, if the number of variables $J$ is large, computational issues can arise while considering all possible pairs of variables. It  restricts the use of their methods to data sets composed of few variables. Moreover, the nature of the approach makes impossible a variable selection.

Selecting variables is challenging in clustering because the role of a variable (relevant or irrelevant for clustering) is defined from the partition which is unobserved. Thus, the selection of the variables and the clustering need to be performed simultaneously. Note that selecting the variables in clustering has two strong benefits: it improves the accuracy of the estimators \citep{AzizyanNIPS2013} because it reduces the number of estimators to be considered and it facilitates the interpretation of the different components  as it only has to be made on the subset of discriminative variables. In a non-model based framework, regularization methods can be used to achieve variable selection in clustering \citep{FriedmanJRSSB2004,PanJMLR2007,WittenJASA2010}. Among these approaches, the sparse K-means \citep{WittenJASA2010} is the most popular because it
requires small computational overhead and is able to manage very
high-dimensional datasets.  The approach uses
a lasso-type penalty to select the set of variables which are relevant for the clustering.
The selection of the number of components is a difficult issue since probabilistic tools are not available and its results are sensitive to the structure of the penalty term. The authors proposed to perform model selection with an extension of the gap statistics \citep{tibshirani2001estimating}. In model-based clustering, \citet{TadesseJASA2005} define  a variable as irrelevant for clustering if its distribution is equal among all the mixture components. A third type of variable, the redundant ones, has been introduced by \citet{RafteryJASA2006}. A redundant variable is a variable that as the same conditional distribution given the relevant variables for each component. However, this type of variable cannot be considered for model \eqref{eq:mixturegrl}, as it requires to model intra-component dependencies. The model-based framework implies that the selection of the variables falls in the scope of model selection. Thus, in the paper, we define the objective of full model selection by the double objective of estimating the number of components $K$ and the subset of relevant variables $\Omega$ as well. In a parametric framework,  variable selection can be performed via information criterion  (see \citet{TadesseJASA2005}, \citet{DeanAISM2010} and \citet{MarbacStatCo2017}) but it leads to computational issues because the number of competing models is of order $2^J$. To circumvent this issue,  \citet{MarbacJoC2019} present a modified EM algorithm \citep{DempsterJRSSB1977,GreenJRSSB1990} that permits to simultaneously perform the variable selection via the Bayesian Information Criterion (BIC; \citet{SchwarzAOS1978}) and the maximum likelihood estimation, for a fixed number of components. Thus, the algorithm only needs to be run for all the possible  numbers of components. Note that all existing methods of variables selection in model-based clustering are restricted to parametric distributions. Thus, if the parametric assumptions are violated, bias can occur for the estimator of the number of components or on the subset of discriminative variables.

This paper addresses the issue of full model selection for non-parametric mixture models defined by \eqref{eq:mixturegrl}. To the best of our knowledge, this paper presents the first method that permits a full-model selection (\emph{i.e.,} estimation of $K$ and $\Omega$) for non-parametric multivariate mixture models. Moreover, it allows many variables to be managed, which makes it  a complementary work to  \citet{KasaharaJRSSB2014} and \citet{Kwon2019estimation}, even in the case where all the variables are considered to be relevant and only the number of components needs to be estimated. As proposed by \citet{TadesseJASA2005}, we consider two types of variables: the relevant variables and the irrelevant variables for clustering. Thus, variable $j$ is said to be irrelevant for clustering if $\eta_{1j}=\ldots=\eta_{Kj}$ and the variable $j$ is said to be relevant for clustering otherwise. A model $M=\{K,\Omega\}$ is defined by  the number of components $K$ and the indices of the relevant variables $\Omega\subset\{1,\ldots,J \}$. Therefore, considering the task of full model selection in \eqref{eq:mixturegrl} implies that each $\bX_i$ is identically distributed according to a non-parametric mixture of $K$ components defined by the density
\begin{equation}\label{eq:modelavecselec1}
g_{M,\bpsi}(\bx_i)=\left(\prod_{j\in\bar\Omega}\eta_{1j}(x_{ij})\right) \left(\sum_{k=1}^K \pi_k \prod_{j \in \Omega} \eta_{kj}(x_{ij}) \right),
\end{equation}
where $\bar\Omega=\{1,\ldots,J\}\setminus\Omega$ contains the indices of the irrelevant variables for clustering and $\bpsi\in\Psi_M$ groups the finite dimensional parameters $\bpi=(\pi_1,\dots,\pi_K)^\top\in\simplex_K$ and the infinite dimensional parameters composed of the univariate densities $\eta_{kj}$.
To achieve the full model selection, we use a discretization of each continuous variable into $B$ bins. The number of bins tends to infinity with the sample size  to ensure the consistency of the approach. Indeed, if $B$ were fixed, the estimated model could be a sub-model of the true model (\emph{i.e.,} the number of components and the subset of the discriminative variables could be underestimated). The distribution of the resulting discretized variables follows a latent class model where each component is a product of multinomial distributions \citep{GoodmanBiometrika1974}. 
This discretization is convenient, because model selection can then be achieved,  for the latent class model, by using the penalized likelihood  (\emph{e.g.,} BIC) whose consistency has been proven for mixture models \citep{KeribinSan2000}. Moreover, in this framework, a specific EM algorithm optimizing the penalized likelihood can be used for simultaneously detecting the subset of the  relevant variables and estimating the model parameters \citep{MarbacJoC2019}, for a known number of components. Thus, by considering an upper-bound of the number of components, the full-model selection can be achieved.
Unlike in \citet{KasaharaJRSSB2014}, the procedure provides a consistent estimation of  the univariate densities of the components $\eta_{kj}$ from the discretized data. Therefore, we prove the consistency of the procedure for a wide range of number of bins $B$, at an appropriate rate that we detail in the paper. 
The consistency of the procedure cannot be proven by using the consistency of information criterion for parametric mixture models \citep{KeribinSan2000} because the parameters space depends on $B$ and thus increases with the sample size. 
The growth rate of $B$ is mainly driven to avoid underestimation of the model while the range of the penalty is mainly driven to avoid overestimation of the model. The case of model underestimation is analyzed by extending the proof of \citet{KeribinSan2000} in order to deal with the increasing dimension of the parameters space. In the case of model overestimation, the asymptotic distribution of the likelihood ratio is investigated by performing a locally conic parametrization \citep{DacunhaESAIM1997,DacunhaAOAS1999} of the model obtained on the discretized data. An upper bound of the likelihood ratio is obtained by controlling, on the one hand, the deviation of the likelihood ratio from its asymptotic distribution by using results on empirical processes stated in \citet{ChernozhukovAOS2014} and, on the other hand, the supremum of the asymptotic distribution by applying deviation results on Gaussian processes \citep{dudley2014uniform}.

The proposed  method uses a discretization that provides an estimator of the densities of the components. However, we advise to use the proposed approach only for model estimation. When the model has been selected, we suggest to use a kernel-based method for density estimation. Indeed, the bin-density estimates are known to be outperformed by kernel-based estimators. Thus, for a real data analysis, we advise to use the proposed approach for model selection then, for the selected model, to perform density estimation with a EM-like algorithm \citep{BenagliaJCGS2009} or by maximizing the smoothed log-likelihood \citep{LevineBiometrika2011}. The final partition is thus computed from the model selected by the proposed method and the densities estimated via a kernel method.

The papers is organized as follows. 
Section~\ref{sec:discret} details the discretization step. 
Section~\ref{sec:converge} states the consistency of the procedure. 
Section~\ref{sec:estimation} presents the algorithm used for the full model selection. 
Section~\ref{sec:num} starts by numerical experiments that compare the proposed approach to standard approaches of variable selection in clustering, then it presents the analysis of benchmark data which illustrates  the relevance of the procedure and introduces some extensions of the approach. 
Section~\ref{sec:concl} gives a conclusion. Mathematical details and supplementary numerical experiments are presented in Supplementary Materials.

\section{Model selection by bin estimation and penalized log-likelihood} \label{sec:discret}
This section considers the estimation of the number of components for model \eqref{eq:modelavecselec1}, from an $n$-sample $\bX_1,\ldots,\bX_n$ with $\bX_i\in\mathcal{X}$ with $\mathcal{X}=\mathcal{X}_1\times\ldots\times\mathcal{X}_J$, $J$ being fixed and $\mathcal{X}_j$ being compact. The method used for selecting the number of components discretizes each variable into $B$ non-overlapping  bins $I_{Bj1},\ldots,I_{BjB}$ such that $\cup_{b=1}^B I_{Bjb}=\mathcal{X}_j$ and for any $(b,b')$ with $b\neq b'$, $I_{Bjb}\cap I_{Bjb'}=\emptyset$. Thus, we consider the function $\sigma_{Bjb}$ with $b\in\{1,\ldots,B\}$, such that $\sigma_{Bjb}(x_{ij})=1$ if $x_{ij}\in I_{jb}$ and $\sigma_{Bjb}(x_{ij})=0$ if $x_{ij}\notin I_{Bjb}$, and we denote by $l_{Bjb}$ the size of the bin $I_{Bjb}$. The discretized variables follow a latent class model where each component is a product of $J$ multinomial distributions each having $B$ levels. Therefore, the pdf of the discretized subject $i$ is
\begin{equation}
f_{M,B,\theta}(\bx_i)=\prod_{j\in \bar\Omega}\prod_{b=1}^B \left(\frac{\alpha_{B1jb}}{l_{Bjb}}\right)^{\sigma_{Bjb}(x_{ij})} \left( \sum_{k=1}^K \pi_k \prod_{j\in \Omega}\prod_{b=1}^B \left(\frac{\alpha_{Bkjb}}{l_{Bjb}}\right)^{\sigma_{Bjb}(x_{ij})}\right),\label{eq:modeldiscret}
\end{equation}
 where $\theta$ groups the component proportions $\pi_k$ and the probabilities $\alpha_{Bkjb}$  that one subject arisen from component $k$ takes level $b$ for the variable $j$ when this variable is discretized into $B$ bins. The parameter space is given by the product of simplexes $S_K \times S_B^{K|\Omega| + (J-|\Omega|)}$, where $|\Omega|$ denotes the cardinal of the set of discriminative variables $\Omega$. Note that $f_{M,B,\theta}$ is an approximation of $g_{M,\bpsi}$ and that this approximation becomes more accurate when $B$ tends to infinity.  The discretized version $f_{M_0,B, \theta_0}$ of the true density $g_0$ will be denoted by $f_{0,B}$.

The probabilities $\alpha_{Bkjb}$ are unknown and must be estimated from the observed sample. This estimation can be achieved by maximizing the log-likelihood defined by
$$
\ell_{n}(f_{M,B,\theta}) = \sum_{i=1}^n \ln f_{M,B,\theta}(\bx_i).
$$
The maximum likelihood statistics for a model with $K$ components and $B$ bins per variables is
$$
T_{n,M,B} = \sup_{\theta \in \Theta_{M,B, \varepsilon}} \ell_{n}(f_{M,B,\theta}),
$$
where, in order to avoid numerical issues, we introduced a threshold $\varepsilon$ such that the parameter space becomes $\Theta_{M,B,\varepsilon}=S_{K,\varepsilon} \times S_{B,\varepsilon}^{K|\Omega| + (J-|\Omega|)}$, with $\varepsilon>0$ being the minimal value of all the elements defined in the simplexes, $\textit{i.e.} \: S_{B,\varepsilon}=\{\boldsymbol{u}\in\mathbb{R}^B: \; u_b >\varepsilon, \; \sum_{b=1}^B u_b=1\}$.   
 Under the condition that $B \varepsilon$ tends to zero as $B$ goes to infinity and $\varepsilon$ to zero, the parameter space $\Theta_{M,B,\varepsilon}$ converges to the whole parameter space. Note that, due to the growth rate of $B$ which will be stated by Assumption~\ref{ass:intervals}\ref{ass:vitesseB} in the next section, it is sufficient to set $\varepsilon^{-1}=O(n^{\alpha + 1})$ for some $\alpha>0$.
This maximization can be achieved via an Expectation-Maximization algorithm (EM algorithm; \citet{DempsterJRSSB1977}). 

The penalized likelihood is defined by substracting from the maximum likelihood statistics a penality term $a_{n,M,B}$ which takes into account the sample size and the complexity of model $M$. Thus, we obtain the following information criterion
\begin{equation}
W_{n,M,B} = T_{n,M,B} - a_{n,M,B}. \label{eq:criterion}
\end{equation}
Depending on the choice of $a_{n,M,B}$ in \eqref{eq:criterion}, different well-known criteria can be considered. Among them one can cite the Akaike criterion (AIC; \citet{AkaikeAISM1970})  or the Bayesian Information Criterion (BIC; \citet{SchwarzAOS1978}) which are  obtained with  $a_{n,M,B}=\nu$ and $a_{n,M,B}=\nu \log(n)/2$ respectively, where $\nu=(K-1) + KJ(B-1)$ is the model complexity.

To select the number of components, we consider the set of competing models $\mathcal{M}$ defined by all the mixture models with at most $K_{\max}$ components and at least three relevant variables (for identifiability reasons), so that
$$
\mathcal{M}=\{M=\{K,\Omega\}: K\leq K_{\max}, \Omega \subset \{1,\ldots,J\}\text{ and } |\Omega| \geq 3\}.
$$
The estimator $\widehat{M}_{n,B}$ of the number of components maximizes the penalized likelihood as follows
$$
\widehat{M}_{n,B} = \argmax_{M \in \mathcal{M}} W_{n,M,B}.
$$
The study of the asymptotic properties of the estimator $\widehat{M}_{n,B}$ is covered by the approach of \citet{KeribinSan2000} only if the number of intervals $B$ does not increase with the sample size $n$. However, in such a case, due to the discretization, the approach would provide an estimator that converges to a model included into the true model. Indeed, we only obtain a lower bound on the number of components and a subset of the discriminative variables. By increasing the number of intervals with $n$, we avoid the issues due to the loss of identifiability. However, we need to investigate the behavior of the statistics $T_{n,M,B}$ and to study  the convergence of $T_{n,M,B}/n$ to the minimum Kullback divergence, which requires controlling empirical processes defined on space having increasing dimension. The next section presents statistical guarantees of the proposed approach.


\section{Convergence in probability of the estimator} \label{sec:converge}

This section investigates the convergence in probability of $\widehat{M}_{n,B}$. It starts by presenting the assumptions required to obtain this convergence, which is then stated. 

\subsection{Assumptions}
The consistency of the estimator is established under four sets of assumptions described below. 
Assumption~\ref{ass:identifiability} and Assumption~\ref{ass:regularite} state the constraints on the model and on the distribution of the components respectively. 
Assumption~\ref{ass:penalty} gives some conditions on the penalty term. Finally, Assumption~\ref{ass:intervals} gives some conditions on the discretization. 

\begin{ass} \label{ass:identifiability}
The number of variables is at least three (\emph{i.e.,} $3\leq J$) and each proportion $\pi_k>0$ is not zero. Moreover, there exists $\Upsilon\subseteq \{1,\ldots,J\}$ such that $|\Upsilon|=3$ and for any $j\in\Upsilon$ the  univariate densities $\eta_{kj}$ are linearly independent.
\end{ass}
\begin{ass}\label{ass:regularite}
\begin{enumerate}[label=(\roman*)]
\item \label{ass:regularite1} There exists a function $\tau$ in $L_1(g_0 \nu)$ such that: $\forall M \in \mathcal{M}$ and $\forall \bpsi \in \Psi_M$, $|\ln g_{M,\bpsi} |<\tau$ $\nu$-a.e.
\item \label{ass:regularite2} There exists a positive constant $L<\infty$ such that $\forall j\in\{1,\ldots ,J \}$ and $\forall x_j\in\mathcal{X}_j$, $|\eta_{kj}'(x_j)|\leq L$. 
\item \label{ass:regularite3} Each variable $j$ is defined on a compact space $\mathcal{X}_j$ and its densities for each component $k$, denoted by $\eta_{kj}$, are striclty positive except on a set of lebesgue measure zero.
\end{enumerate}
\end{ass}
\begin{ass}\label{ass:penalty}
\begin{enumerate}[label=(\roman*)]
\item \label{ass:penalty1} $a_{n,M,B}$ is an increasing function of $K$, $|\Omega|$ and $B$.
\item \label{ass:penalty2} For any model $M$, $a_{n,M,B}/n$ tends to 0 as $n$ tends to infinity.
\item \label{ass:penalty3} For any model $M$, $B/a_{n,M,B}$ tends to 0 as $n$ tends to infinity.
\item \label{ass:penalty4} For any models $M$ and $\widetilde{M}$ with $M\subset \widetilde{M}$, $a_{n,\widetilde M,B}/a_{n,M,B}$ tends to infinity as $n$ tends to infinity.
\end{enumerate} 
\end{ass}
\begin{ass}\label{ass:intervals}
\begin{enumerate}[label=(\roman*)]
\item The number of bins $B$ tends to infinity with $n$ in the following way $\lim_{n\to\infty} B=\infty$ and $\lim_{n\to\infty} B(\ln^3 n)/n=0$.\label{ass:vitesseB}
\item \label{ass:intervals2} 
The length of the each interval is not zero and  satisfies, for all $j\in\{1,\ldots ,J \}$ and $b\in\{1,\ldots ,B \}$,  $l_{Bjb}^{-1}=O(B)$.
\item \label{ass:intervals1} Let $\mathcal{I}_{jB}$ be the set of the upper bounds of the $B$ intervals, then, for any value $x_j\in\mathcal{X}_j$,  $d(x_j,\mathcal{I}_{jB})$ tends to zero as $B$ tends to infinity. 
\end{enumerate}
\end{ass}

Assumption~\ref{ass:identifiability} is derived from the conditions of identifiability for finite mixtures of nonparametric measure products (see Theorems 8 and 9 in \citet{AllmanAOS2009}). Because Theorems 8 and 9 in \citet{AllmanAOS2009} consider all the variables as revelant for clustering, we need to extend their assumptions such that there are at least three revelant variables to obtain the identifiability of the model \eqref{eq:modelavecselec1}.

Assumption~\ref{ass:regularite} gives sufficient conditions on the component distributions to ensure that the results of \citet{DacunhaAOAS1999} can be applied to the mixture model obtained after discretization. 

Assumption~\ref{ass:penalty} presents standard conditions for penalized likelihood model selection in the case of embedding models. It generalizes the usual conditions for selecting the number of components \citep{KeribinSan2000,GassiatAIHP2002,ChambazAOS2006} to the case of feature selection for mixture models. Conditions \ref{ass:penalty1} and \ref{ass:penalty3} permit avoiding the overestimation of the model (\emph{i.e.,} overestimation of the number of components or of the support of the relevant variables), while condition  \ref{ass:penalty2} permits tavoiding the underestimation of the model by making the penalty term negligeable with respect to the model bias. Note that Assumption~\ref{ass:penalty}  allows the BIC penalty to be considered.

Even if Assumption~\ref{ass:identifiability}  provides  the identifiability of model \eqref{eq:modelavecselec1}, after the discretization, model \eqref{eq:modeldiscret} could be not identifiable if the number of intervals $B$ were fixed. 
As an example, one can consider a bi-component mixture model with equal proportions  defined with a first component following a product of $J\geq 3$ beta distributions $\mathcal{B}e(\alpha,\alpha)$ and a second component following a product of $J\geq 3$  beta distributions $\mathcal{B}e(2\alpha,2\alpha)$, with $\alpha\geq 1$. This model is identifiable but the model \eqref{eq:modeldiscret} defined after the discretization of each variable into two bins of equal size (\emph{e.g.,} for any $j$, $\sigma_{j1}(u)=1$ if $0\leq u\leq 1/2$, $\sigma_{j1}(u)=0$ if $1/2<u\leq 1$ and $\sigma_{j2}(u)=1-\sigma_{j1}$) is not identifiable (\emph{i.e.,} the two mixture components follow the same distribution for the discretized data). However, if the number of bins is strictly larger than two and if each interval has a length not equal to zero, then the model \eqref{eq:modeldiscret} becomes identifiable.

The model identifiability is obtained by Assumption~\ref{ass:intervals} that states conditions on the discretization. In particular,  the number of levels has to tend to infinity when the sample size increases in such a way that the size of the largest interval tends to zero when the sample size tends to infinity. However, its growth rate has to be  upper bounded, which is a key point to control the convergence of the estimators. Note that Assumption~\ref{ass:intervals}\ref{ass:intervals1} uses the same ideas as Lemma 17 in \citet{AllmanAOS2009} and that this condition is not stringent. For instance the bounds of the intervals can be determined by the quantiles $1/B,\ldots,B/B$.  In addition, the sizes $l_{jb}$ can vary from one bin to another. This is for instance the case when we consider the quantiles.  However, we cannot allow a bin to be exponentially small with $n$, in order to keep the asymptotic behavior of our estimator which is stated in the next subsection.  Note that Assumption~\ref{ass:intervals} allows to consider the rate $B=n^{1/3}$ that is usual for bin-density estimation.

Finally, the assumption on the compactness of $\mathcal{X}_j$ can be relaxed if some densities defined on $\mathbb{R}$ are wanted to be considered. In such case, the estimates of the densities are considered on the compact $[\min_i x_{ij},\max_i x_{ij}]$ defined from the observed sample, and the estimates of the densities are zero outside this interval.

\subsection{Convergence in probability of the estimator}
We state the consistency of the estimator $\widehat{M}_{n,B}$ then we give its (sketch of) proof. Note that the proof of all the numbered equations are given in Section~\ref{sec:maths} of the Supplementary Materials. Finally, we explain the key points of the proof which are different from the proof of the consistency of information criteria for parametric mixture models stated in \citet{KeribinSan2000}.

\begin{theorem} \label{thm:cvproba}
Assume that independent data arise from \eqref{eq:modelavecselec1} with the true model $M_0=\{K_0,\Omega_0\}$, that Assumptions~\ref{ass:identifiability}, \ref{ass:regularite}, \ref{ass:penalty} and \ref{ass:intervals} hold true, and that the set of competing models $\mathcal{M}$ is defined with a known upper bound for the number of clusters $K_{\max}$. Then, $\widehat{M}_{n,B}$ converges in probability to $M_0$.
\end{theorem}

\begin{proof}[Proof of Theorem~\ref{thm:cvproba}]
The proof is divided into three parts: the case where $M_0$ is underestimated (\emph{i.e.,} $K<K_0$ or $\Omega_0 \not\subseteq \Omega$), the case where the subset of the relevant variables is overestimated with $K_0$ (\emph{i.e.,} $K=K_0$ and $\Omega_0\varsubsetneq\Omega$) and the case where the number of components and the subset of relevant variables are overestimated (\emph{i.e.,} $K>K_0$ and $\Omega_0\subseteq \Omega$).\\
$\bullet$ \underline{Part 1:} We consider the case where   $M_0$ is underestimated. Thus, we consider the set of  models $$\mathcal{N}_1=\{M=\{K,\Omega\}\in \mathcal{M}: K\leq K_0 \text{ or } \Omega_0 \not\subseteq \Omega\}.$$ 
The probability to underestimate the true model can be upperbounded as follows
$$\mathbb{P}(\widehat{M}_{n,B} \in \mathcal{N}_1) \leq \sum_{M \in \mathcal{N}_1} \mathbb{P}(W_{n,M,B}-W_{n,M_0,B} \geq 0). $$
For any $g_{M,\bpsi}$ given by model \eqref{eq:modelavecselec1},  Assumption~\ref{ass:regularite}\ref{ass:regularite1} implies that $\mathbb{E}_{g_0}[\ln g_{M,\bpsi}]$ is defined. The Kullback-Leibler divergence from model $M$ to the true distribution $g_0$ is  defined by
$$
\text{KL}(g_0,\mathcal{G}_M): =\inf_{\bpsi \in \Psi_M}\mathbb{E}_{g_0}\left[\ln \frac{g_0}{g_{M,\bpsi}}\right].
$$
Using the definition of $\mathcal{N}_1$ and the identifiability of $g_0$ (ensured by Assumption~\ref{ass:identifiability}),  for each $M\in\mathcal{N}_1$, there exists some $\delta_M>0$ such that $\text{KL}(g_0,\mathcal{G}_M)\geq\delta_M$. In Section~\ref{proof:cvprobaunderestimation} of the Supplementary Materials, we prove the following asymptotic bound in probability:
\begin{equation} \label{eq:cvprobaunderestimation}
\frac{1}{n} \left(T_{n,M,B} - \ell_n(g_0) \right) \leq  -\delta_M + o_{\mathbb{P}}(1).
\end{equation}
This combined to the properties of the penalty (see Assumption~\ref{ass:penalty}) implies that for any $M\in\mathcal{N}_1$
$$
\frac{1}{n}\left(W_{n,M,B} - W_{n,M_0,B}\right) \leq -\delta_M + o_\mathbb{P}(1).
$$
Therefore, noting that $\delta_M>0$ and that the cardinal of $\mathcal{N}_1$ is fixed and finite, we have
\begin{equation} \label{eq:nonunder}
\lim_{n\to\infty}  \mathbb{P}(\widehat{M}_{n,B} \in \mathcal{N}_1) = 0.
\end{equation}
Thus, the probability of underestimating the model tends to zero as $n$ tends to infinity.\\
$\bullet$ \underline{Part 2:} We consider the case where the number of components is correct but the subset of the relevant variables is overestimated. Thus, we consider a model  $M\in \mathcal{N}_2$ where $$\mathcal{N}_2=\{M=\{K,\Omega\}\in \mathcal{M} : K = K_0 \text{ and } \Omega_0 \varsubsetneq \Omega\}.$$ 
We have the following upper-bound
$$
\mathbb{P}(\widehat{M} \in \mathcal{N}_2) \leq \sum_{M \in \mathcal{N}_2} \mathbb{P}(W_{n,M,B} \geq W_{n,M_0,B})= \sum_{M \in \mathcal{N}_2} \mathbb{P}\left(\frac{T_{n,M,B} - T_{n,M_0,B}}{a_{n,M_0,B}} \geq \frac{a_{n,M,B}}{a_{n,M_0,B}} - 1\right).
$$
Note that for any $M \in \mathcal{N}_2$, we have $\delta_M=0$ and thus the reasonning used to demonstrate that $M_0$ is not underestimated cannot be used. However, the models on $\mathcal{N}_2$ are identifiable because $K=K_0$. Thus, using usual results on likelihood ratio, for a fixed value of $B$, $ 2 \left( T_{n,M,B} -  T_{n,M_0,B} \right)$ is asymptotically distributed like a $\chi^2(\Delta) $ where $\Delta$ is given by the difference of the dimensions of both parameter spaces $\Delta = (B-1) (K-1) (|\Omega|-|\Omega_0|)$. As $\Delta$ goes to infinity with $B$, thus with $n$, we have the following asymptotic distribution
$$\frac{1}{\sqrt{2\Delta}} \left[2 \left( T_{n,M,B} -  T_{n,M_0,B} \right) - \Delta \right] \xrightarrow{d} \mathcal{N}(0,1).$$
We rewrite
$$\frac{T_{n,M,B} - T_{n,M_0,B}}{a_{n,M_0,B}} = \frac{1}{a_{n,M_0,B}} \sqrt{\frac{\Delta}{2}} \left( \frac{1}{\sqrt{2\Delta}} \left[2 \left( T_{n,M,B} -  T_{n,M_0,B} \right) - \Delta \right]\right) + \frac{\Delta}{2 a_{n,M_0,B}},$$
and conclude, by making use of Slutsky's lemma and Assumption~\ref{ass:penalty}~\ref{ass:penalty3}, that
$$\frac{T_{n,M,B} - T_{n,M_0,B}}{a_{n,M_0,B}} = o_{\mathbb{P}}(1).$$
For any $M\in\mathcal{N}_2$, Assumption~\ref{ass:penalty}~\ref{ass:penalty4} implies that $ a_{n,M,B}/a_{n,M_0,B} - 1>0$, thus, as the cardinal of $\mathcal{N}_2$ is finite and does not depend on $B$, we can conclude that
\begin{equation} \label{eq:probover1}
\mathbb{P}(\widehat{M} \in \mathcal{N}_2) = 0.
\end{equation}
$\bullet$ \underline{Part 3:} We consider the case where the number of components and the subset of the relevant variables are overestimated. Thus, we consider a model  $M\in \mathcal{N}_3$ where
$$
\mathcal{N}_3 =\{M=\{K,\Omega\}\in \mathcal{M}: K > K_0 \text{ and } \Omega_0 \subseteq \Omega\}.
$$
Note that $\mathcal{N}_3 = \mathcal{M}\setminus\{\mathcal{N}_1 \cup \mathcal{N}_2 \cup M_0\}$. 
The probability of overestimating the model (\emph{i.e.,} $\widehat{M} \in \mathcal{N}_3$) can be upperbounded by
$$
\mathbb{P}(\widehat{M} \in \mathcal{N}_3) \leq \sum_{M \in \mathcal{N}_3} \mathbb{P}(W_{n,M,B} \geq W_{n,M_0,B}).
$$
Note that for any $M \in \mathcal{N}_3$, we have $\delta_M=0$ and thus the reasonning used to demonstrate that $M_0$ is not underestimated cannot be used. Moreover, because for $M \in \mathcal{N}_3$, $K>K_0$, the model suffers from a loss of identifiability implying that the likelihood ratio must be carrefully investigated. We have, for any $M\in\mathcal{N}_3$
$$
\mathbb{P}\left(W_{n,M,B} \geq W_{n,M_0,B}\right) = \mathbb{P}\left(\frac{T_{n,M,B} - T_{n,M_0,B}}{a_{n,M_0,B}} \geq \frac{a_{n,M,B}}{a_{n,M_0,B}} - 1\right).
$$
Applying the locally-conic parametrization proposed by \citet{DacunhaESAIM1997,DacunhaAOAS1999} on model \eqref{eq:modeldiscret}, and noting that Assumption~\ref{ass:regularite} holds true, we can rewrite the log-likelihood ratio as in the proof of Lemma 3.3 in \cite{KeribinSan2000}
\begin{multline}\label{eq:cvover}
T_{n,M,B} - \ell_n f_{0,B} \\ = \sup \left\lbrace \sup_{d \in D_B} \frac{1}{2} \mathcal{G}_n^2(d) \mathds{1}_{\mathcal{G}_n(d)\geq 0} ;
\sup_{d_1 \in D_{1B}, d_2 \in D_{2B}} \frac{1}{2} \left(\mathcal{G}_n^2(d_1) +  \mathcal{G}_n^2(d_2) \mathds{1}_{\mathcal{G}_n(d_2)\geq 0} \right) \right\rbrace  (1 + o_{\mathbb{P}}(1))
\end{multline}
where, for each function $d$, $\mathcal{G}_n(d)=n^{-1/2} \sum_{i=1}^n d(X_i)$; the considered spaces of functions as well as the definition of $f_{0,B}$ are detailled in Section~\ref{sec:cvover} of the Supplementary Material.
Note that 
$$\sup_{d \in D_B} \frac{1}{2} \mathcal{G}_n^2(d) \mathds{1}_{\mathcal{G}_n(d)\geq 0} =  \frac{1}{2} \left(\sup_{d \in D_B} \mathcal{G}_n(d)\right)^2.$$
In addition, as $D_{1B}$ and $D_{2B}$ are subspaces of $D_B$, we have 
$$ \sup_{d_1 \in D_{1B}, d_2 \in D_{2B}} \frac{1}{2} \left(\mathcal{G}_n^2(d) +  \mathcal{G}_n^2(d) \mathds{1}_{\mathcal{G}_n(d)\geq 0} \right) \leq \left(\sup_{d\in D_{B,s}} \mathcal{G}_n(d)\right)^2, $$
where $D_{B,s}$ is the symmetrized space $D_B \cup (-D_B)$. Therefore, we deduce that $$T_{n,M,B} - \ell_n f_{0B} \leq \left(\sup_{d\in D_{B,s}} \mathcal{G}_n(d)\right)^2 ( 1 + o_{\mathbb{P}}(1)).$$
Thus, using the fact that $D_{B,s}$ is a symmetric space, we obtain that, for any $\varepsilon >0$, for $n$ sufficently large,
$$\left\lbrace \frac{T_{n,M,B} - \ell_n f_{0,B}}{a_{n,M_0,B}} > 4 \varepsilon \right\rbrace \subset \left\lbrace \left|\sup_{d\in D_{B,s}} \mathcal{G}_n(d) \right| > 2 \sqrt{\varepsilon a_{n,M_0,B}} \right\rbrace.$$
It implies that,  
$$\mathbb{P}\left(\frac{T_{n,M,B} - \ell_n f_{0,B}}{a_{n,M_0,B}} >\varepsilon\right) \leq \mathbb{P}\left(\sup_{d\in D_{B,s}} \xi_d > \sqrt{\varepsilon a_{n,M_0,B}}\right) 
 + \mathbb{P}\left(\left|\sup_{d\in D_{B,s}} \mathcal{G}_n(d) - \sup_{d\in D_{B,s}} \xi_d \right|>\sqrt{\varepsilon a_{n,M_0,B}}\right),
$$
where $(\xi_d)_d$ is a Gaussian process indexed by $D_{B,s}$, with covariance the usual Hilbertian product on $L^2_{f_{0,B}}$. Note that under our Assumptions, for a fixed value $B^{\star}$ of $B$, $\sup_{d \in D_{B^{\star},s}} \mathcal{G}_n(d) $ converges in distribution to $\sup_{d \in D_{B^{\star},s}} \xi_d$ (see Section~\ref{sec:cvover} in the Supplementary Materials) as $n$ goes to infinity. However, as in our context, $B$ goes to infinity with $n$, we need to control its influence on the deviations.

Under Assumption~\ref{ass:regularite}, we will control the first term on the right-hand side by using existing deviation bounds for the supremum of Gaussian processes (see Section~\ref{sec:asymptV} for details), which will lead to
\begin{equation} \label{eq:asymptV}
\lim_{n\to\infty} \mathbb{P}\left(\sup_{d\in D_{B,s}} \xi_d > \sqrt{\varepsilon a_{n,M_0,B}}\right)  =0.
\end{equation}
Moreover, using the results of the approximation of suprema of general empirical processes by a sequence of suprema of Gaussian processes \citep{ChernozhukovAOS2014}, we obtain (see Section~\ref{sec:asymptconcentre})
\begin{equation}\label{eq:asymptconcentre}
\lim_{n\to\infty} \mathbb{P}\left(\left|\sup_{d\in D_{B,s}} \mathcal{G}_n(d) - \sup_{d\in D_{B,s}} \xi_d \right|>\sqrt{\varepsilon a_{n,M_0,B}}\right) =0.
\end{equation}
Thus, we have for any $\varepsilon>0$, 
\begin{equation*}
\lim_{n\to\infty} \mathbb{P}\left(\frac{T_{n,M,B} - \ell_n f_{0,B}}{a_{n,M_0,B}} >\varepsilon\right)=0.
\end{equation*}
 Noting that for any $M\in\mathcal{N}_3$, Assumption~\ref{ass:penalty} implies that $ a_{n,M,B}/a_{n,M_0,B} - 1>0$ and noting that the cardinal of $\mathcal{N}_3$ is finite and fixed (it does not depend on $B$) we can then conclude that
\begin{equation} \label{eq:probover}
\mathbb{P}(\widehat{M} \in \mathcal{N}_3) = 0.
\end{equation}
Combining equations \eqref{eq:nonunder}, \eqref{eq:probover1} and  \eqref{eq:probover} leads to the convergence in probability of $\widehat{M}$ to $M_0$.
\end{proof}

\paragraph*{Some comments:}
Note that the arguments used in \citet{KeribinSan2000} to prove that underestimation is avoided, cannot be used here. Indeed, the proof of Theorem 2.1 in  \citet{KeribinSan2000} considers parameters that are defined on fixed dimensional space and thus cannot be used to obtain \eqref{eq:cvprobaunderestimation}. In our context, we require that $B$ tends to infinity with $n$ (see Assumption~\ref{ass:intervals}) to ensure the identifiability and thus the convergence of $\inf_{\theta}\mathbb{E}_{g_0}\left[\ln f_{M_0,B,\theta_{M_0,B}^\star} - \ln f_{M,B,\theta}\right]$ to a quantity lower-bounded by $\delta_M$ where $\theta_{M_0,B}^\star=\argmax_{\theta \in \Theta_{M_0,B}}\mathbb{E}_{g_0}\left[\ln f_{M_0,B,\theta}\right]$ (this convergence is ensured by Assumption~\ref{ass:regularite}\ref{ass:regularite2}, as discussed in the proof of \eqref{eq:cvprobaunderestimation}).  Note that the existence of $\theta_{M_0,B}^\star$ is ensured by the fact that $\Theta_{M_0,B}$ is compact and that the Kullback-Leibler divergence is continuous. 

Note also that the arguments used in \citet{KeribinSan2000} to prove that overestimation is avoided cannot be used here either. Indeed, despite the fact that $T_{n,M,B^\star} - \ell_n f_{0B^\star}$ converges in distribution for a fixed $B^\star$ and that $1/a_{n,M,B^\star}$ tends to 0, we cannot directly conclude that  $ \mathbb{P}\left(\left[T_{n,M,B} - T_{n,M_0,B}\right]/a_{n,M,B} > \varepsilon\right)$ tends to zero as $n$ tends to infinity, for any $\varepsilon>0$ and any $M\in\mathcal{N}_3$ because we require that $B$ tends to infinity with $n$ to avoid the underestimation (see Assumption~\ref{ass:intervals}).  

\section{Estimation of the best model} \label{sec:estimation} 

The estimation of $\widehat{M}_{n,B}$ requires an optimization over a discrete space whose cardinal is of order $2^J K_{\max}$. Thus, an exhaustive approach computing $W_{n,M,B}$ for each $M$ in $\mathcal{M}$ is not doable in practice. As the combinatorial issue is mainly due to the feature selection, we follow the approach of \citet{MarbacJoC2019} that consists of simultaneously performing feature selection and parameter estimation, with a fixed number of components, via a specific EM algorithm optimizing the penalized likelihood. Thus, for a fixed value of $K$, the goal of the algorithm is to estimate
$$
\widehat{M}_{n,B,K}=\argmax_{\{M=(K,\Omega) \text{ with } \Omega \subseteq \{1,\ldots,J\} \text{ and }  |\Omega|\geq 3 \}} W_{n,M,B}.
$$
The following EM algorithm permits the estimation of the model parameters and the detection of the subset of relevant variables, for a fixed number of components $K$. 
This algorithm ensures that the penalized log-likelihood increases at each iteration. Thus,
parameter estimation is achieved by maximum likelihood and model selection is done with an information criterion with penalty $a_{n,K,\Omega,M}=\nu_{K,\Omega,M} c_n$ where $\nu_{K,\Omega,M}=(K-1) + |\Omega|K(B-1) + (J-|\Omega|)(B-1)$ is the number of model parameters. The algorithm considers a fixed number of components $K$ and starts at an initial point $\{\Omega^{[0]},\theta^{[0]}\}$. Its iteration $[r]$ is composed of two steps:\\
\textbf{E-step} Computation of the fuzzy partition 
$$t_{ik}^{[r]}:=\dfrac{\pi_k^{[r-1]} \prod_{j \in \Omega^{[r-1]}} \prod_{b=1}^B \left(\alpha_{Bkjb}^{[r-1]}\right)^{\sigma_{Bjb}(x_{ij})} }{\sum_{\ell=1}^K\pi_{\ell}^{[r-1]} \prod_{j \in \Omega^{[r-1]}} \prod_{b=1}^B \left(\alpha_{B\ell jb}^{[r-1]}\right)^{\sigma_{Bjb}(x_{ij})}},$$
\textbf{M-step} Maximization of the expectation of the penalized complete-data log-likelihood over $\Omega$ and $\theta$ such 
$$
\Omega^{[r]}=\{j: \Delta_j^{[r]}>0\},\; \pi_k^{[r]}=\dfrac{n_k^{[r]}}{n} \text{ and } \alpha^{[r]}_{Bkjb}=\left\{ \begin{array}{rl}
\tilde \alpha^{[r]}_{Bkjb} & \text{if } j \in \Omega^{[r]} \\
\bar{\alpha}_{Bkjb} & \text{otherwise}
\end{array}\right.,$$
where $$\Delta_j^{[r]}= \sum_{i=1}^n\sum_{b=1}^B \sigma_{Bjb}(x_{ij}) \sum_{k=1}^K t_{ik}^{[r]} \ln \left(\frac{\tilde \alpha^{[r]}_{Bkjb}}{ \bar{\alpha}_{Bkjb}}\right)- (K-1)(B-1) c_n$$ is the difference between the maximum of the expected value of the penalized complete-data log-likelihood obtained when variable $j$ is relevant and when it is irrelevant, with 
$$\tilde \alpha^{[r]}_{Bkjb}=\frac{1}{n_k^{[r]}} \sum_{i=1}^nt_{ik}^{[r]} \sigma_{jb}(x_{ij}), \; \; \;  \bar{\alpha}_{Bkjb}=\frac{1}{n}\sum_{i=1}^n \sigma_{jb}(x_{ij}) \; \; \; \text{ and } n_k^{[r]}=\sum_{i=1}^n t_{ik}^{[r]}.$$
Note that, when less than three variables happen to have a positive value for $\Delta_j^{[r]}$, the M-step selects in $\Omega^{[r]}$ the three variables having the largest values of $\Delta_j^{[r]}$.
 To obtain the pair $\Omega$ and $\theta$ maximizing the penalized observed-data log-likelihood, for a fixed number of components, many random initializations of this algorithm should be done. Hence, the model (\emph{i.e.,} $K$ and $\Omega$) and the parameters maximizing the penalized observed-data log-likelihood are obtained by performing this algorithm for every values of $K$ between $1$ and $K_{\max}$. By considering $c_n=(\ln n)/2$, this algorithm carries out the model selection according to the BIC.

From the previous algorithm, we obtain an estimator of the model and of its parameters. Indeed, $\hat{\alpha}_{kjb}/l_{jb}$ estimates the density $\eta_{kj}(u)$ for any $u$ such that $\sigma_{jb}(u)=1$. However, the bin-based density estimators are generally outperformed by kernel-based estimators. Thus, we advise to use the proposed approach only for model estimation. Then, for the selected model, kernel-based density estimates provided by the EM-like algorithm \citep{BenagliaJCGS2009} or by maximizing the smoothed log-likelihood \citep{LevineBiometrika2011} should be considered. However, note that establishing asymptotic properties of those kernel-based density estimators is still an open question.

\section{Numerical experiments} \label{sec:num}

\subsection{Comparing the methods for selecting the number of components}
In this section, we assess the performance of our estimator of the number of components without considering the task of variable selection (\emph{i.e.,} $\Omega=\{1,\ldots,J\}$ is supposed to be known). We consider the estimator $\hat K$ obtained with a BIC penalty and $B$ levels defined by the empirical quantiles $1/B$,...,$B/B$ implying that  $\sigma_{jb}(u)=1$ only if $u\in[q_{jb-1}, q_{jb}]$ and $\sigma_{jb}(u)=0$ otherwise, with $q_{j0}=\min_i x_{ij}$, $q_{jB}=\max_{i} x_{ij}$ and $q_{jb}$ the empirical quantile of order $b/B$ for $b=1,\ldots,B-1$. Three values of $B$ are investigated: $[n^{1/5}]$, $[n^{1/6}]$ and $[n^{1/7}]$, where $[.]$ denotes the closest integer. Our estimator is compared to the procedure of \citet{Kwon2019estimation}: SVT and to the procedures of \citet{KasaharaJRSSB2014}: max ave-rk$^+$ statistics with $M_0=4$. For the comparison, we consider the simulation setup of \citet{KasaharaJRSSB2014}. This setup originally considers four designs with two variables and one design with eight variables (Design 5). Our approach cannot be used when only two variables are observed for identifiability issues. This is a drawback of the proposed method that is due to the discretization step. 
Indeed, \citet{KasaharaJRSSB2014} show that the mixture of products of two densities is identifiable under mild assumption. However, this is no longer the case for a mixture of products of two multinomial distributions. Note that the approaches of \citet{KasaharaJRSSB2014} and  \citet{Kwon2019estimation} are based on estimator of bivariate density and we advise to use these methods for such data (in this case variable selection does not make sense). Hence, these methods suffer from computational issue if the number of variables is large (explaining that we do not run these methods on the next section) while the proposed method does not suffer from this problem. This illustrates the complementarity of the three approaches. We simulate 1000 samples of size $n=500$ and of size $n=2000$ from a mixture of three Gaussian distributions (see numerical experiments in \citet{KasaharaJRSSB2014} and  \citet{Kwon2019estimation}) defined with equal proportions, centers $\mu_1=(0,0,0,0,0,0,0,0)^\top$, $\mu_2=(1,2,0.5,1,0.75,1.25,0.25,0.5)^\top$ and $\mu_3=(2,1,1,0.5,1.25,0.75,0.5,0.25)^\top$ and covariance matrices equal to the identity matrix of size eight.

\begin{table}[ht!]
\centering
\begin{tabular}{ccccccccc}
  \hline
  & \multicolumn{4}{c}{$n=500$}
  & \multicolumn{4}{c}{$n=2000$}\\
 & $K=1$ & $K=2$ & $K=3$ & $K\geq 4$  & $K=1$ & $K=2$ & $K=3$ & $K\geq 4$ \\ 
  \hline
BIC with $B=n^{1/7}$ & 0.000 & 0.996 & 0.004 & 0.000 & 0.000 & 0.517 & 0.483 & 0.000 \\ 
  BIC with $B=n^{1/6}$ & 0.001 & 0.999 & 0.000 & 0.000 & 0.000 & 0.876 & 0.124 & 0.000 \\ 
  BIC with $B=n^{1/5}$ & 0.001 & 0.999 & 0.000 & 0.000 & 0.000 & 0.992 & 0.008 & 0.000 \\ 
  \hline
  SVT & 0.000 & 0.992 & 0.008 & 0.000 & 0.000 & 0.493 & 0.507 & 0.000 \\ 
  ave-rk & 0.142 & 0.810 & 0.047 & 0.001 & 0.005 & 0.776 & 0.214 & 0.005 \\ 
  AIC by ave-rk & 0.012 & 0.867 & 0.119 & 0.003 & 0.000 & 0.587 & 0.399 & 0.013 \\ 
  BIC by ave-rk & 0.284 & 0.715 & 0.001 & 0.000 & 0.035 & 0.942 & 0.023 & 0.000 \\ 
  HQ by ave-rk & 0.078 & 0.909 & 0.013 & 0.000 & 0.004 & 0.878 & 0.117 & 0.001 \\  
   \hline
\end{tabular}
    \caption{Empirical probabilities of selecting the different numbers of clusters.} \label{tabcomparing}
\end{table}

Table~\ref{tabcomparing} presents the outcome of the simulation. Note that the results of the approaches of \citet{KasaharaJRSSB2014} and  \citet{Kwon2019estimation} arise from Section 4 in \citet{Kwon2019estimation} and thus are not obtained on the same samples.  For $n=500$, all the procedures performs poorly by mainly selecting two components. The best approach is the AIC (ave-risk) procedure. The proposed approach obtains similar results than the SVT approach which are worse than those of the AIC procedure. However, when the sample size increases ($n=2000$), the SVT approach obtains the best results (the true number of components is detected 50.7$\%$ of time). The proposed method obtains also good results (the true number of components is detected 48.3$\%$ of time) when $B=[n^{1/7}]$ but its results deteriorate when the number of levels grows faster (\emph{i.e.,} $B=[n^{1/6}]$ and $B=[n^{1/5}]$). This illustrates that the results of the proposed method are sensitive to the choice of the number of levels used for the discretization, for a fix sample size. However, all the growth rates of the number of levels respecting Assumption~\ref{ass:intervals}.\ref{ass:vitesseB} provide consistent estimator as illustrated by the experiments presented in Section~\ref{sec:extraexpe} of the Supplementary Materials.

\subsection{Comparing the methods for a full model selection}  \label{sec:simucomparing}

This section compares approaches for a full model selection (\emph{i.e.,} estimation of the subset of the relevant variables and on the number of components) on simulated data. We compare the proposed approach  with a BIC applied on a Gaussian mixture model, with the sparse $K$-means approach  and with the non-parametric approach considering all the variables as relevant. 
The results of the proposed clustering method are obtained by performing full model selection with $B$ levels defined by the empirical quantiles $1/B$,...,$B/B$ where $B=[n^{1/6}]$ and a BIC like penalty and then by estimating  the mixture components for the selected model by maximizing the smoothed log-likelihood with a bandwidth, for variable $j$, equal to $\hat\sigma_j n^{-1/5}$ where $\hat\sigma_j$ is the empirical standard deviation of variable $j$. Thus, when the discretization is performed, the model selection can be achieved via the R package VarSelLCM \citep{MarbacBioinfo2019} then, when the best model is selected, the maximization of the smoothed log-likelihood is achieved via the R package mixtools \citep{BenagliaJSS2009}. The parametric mixture model considers that all the components are Gaussian (this approach is also implemented in the R package VarSelLCM) and uses the BIC to perform model selection. The sparse $K$-means approach is implemented in the R package sparcl \citep{WittenJASA2010} and consists in the sparse $K$-means algorithm initialized with the partition provided by the sparse hierchical ascendant classification with the ``average" method.
 The sparse $K$-means estimates weights for the variables and thus we consider that a variable is estimated as relevant if its weight is more than $\iota$ and that a variable is irrelevant if its weight is less than $\iota$, where  the small threshold $\iota=10^{-6}$ is introduced to avoid numerical issues.  
 Finally, the non-parametric mixture model is implemented the R package mixtools \citep{BenagliaJSS2009} and considers the estimator maximizing the smoothed log-likelihood with a bandwidth, for variable $j$, equal to $\hat\sigma_j n^{-1/5}$ where $\hat\sigma_j$ is the empirical standard deviation of variable $j$.
To compare the different methods of clustering, we generate data from a mixture with three components and equal proportions ($\pi_k=1/3$). The density of $X_i$ given $Z_i$ is a product of univariate densities such that $X_{ij}=  \sum_{k=1}^K z_{ik}\delta_{kj} + \xi_{ij}$ where all the $\xi_{ij}$ are independent and where  $\delta_{11}=\delta_{12}=\delta_{23}=\delta_{24}=\delta_{35}=\delta_{36}=\tau$, while all remaining $\delta_{kj}=0$,  which implies that only the first six variables are relevant for clustering. Three distributions are considered for the $\xi_{ij}$ (standard Gaussian, Student with three degrees of freedom and Laplace) and the value of $\tau$ is defined to obtain a theoretical misclassification rate of $5\%$ ($\tau$ is equal to 1.94, 2.60 and 2.52 for the Gaussian, Student and Laplace distributions respectively). In the Section~\ref{sec:extraexpe} of the Supplementary Materials., all the experiments are also run with theoretical  misclassification rates equal to $10\%$ and $15\%$.



 

\paragraph{Selection of the discriminative features}
To investigate the performances of the competing methods for feature selection, we first consider the situation with a known number of components. Thus, the model selection consists in performing the feature selection.  We consider the methods which automatically provide an estimator of the relevant variables (\emph{i.e.,} the proposed method, sparse $K$-means and VarSelLCM). Accuracy of this selection is measured by sensitivity (probability to detect as relevant a true discriminative variable) and specificity (probability to detect as irrelevant a true non discriminative variable). Table~\ref{tab:simu1feature.sen} and \ref{tab:simu1feature.spe} present the sensibility and the specificity obtained by the proposed approach and the parametric approach. They exhibit an advantage of the parametric method when the distribution is well-specified, but only for small samples ($n=100$). The reason is that, for such samples, the proposed method only finds a part of the relevant variables. However, both methods perform well for larger samples. Moreover, the proposed method obtains similar results for the two other distributions of the components while the results of the parametric approach are strongly deteriorated for both sensibility and specificity, especially for heavy tailed distributions (Student distribution). 

\begin{table}[ht] 
\centering
\begin{tabular}{ccccccccccc}
  \hline 
& & \multicolumn{3}{c}{Proposed method}& \multicolumn{3}{c}{VarSelLCM} & \multicolumn{3}{c}{Sparcl} \\
& & \multicolumn{3}{c}{$n$}& \multicolumn{3}{c}{$n$} & \multicolumn{3}{c}{$n$} \\
Component & $J$ &  100 & 250 & 500& 100 & 250 & 500& 100 & 250 & 500 \\
  \hline
  Gaussian & 20 & 0.81 & 1.00 & 1.00 & 1.00 & 1.00 & 1.00 & 0.86 & 0.91 & 0.95 \\ 
   & 50 & 0.66 & 0.99 & 1.00 & 0.96 & 1.00 & 1.00 & 0.91 & 0.95 & 0.97 \\ 
     & 100 & 0.35 & 0.70 & 1.00 & 0.69 & 1.00 & 1.00 & 0.78 & 0.97 & 0.98 \\ 
  Student & 20 & 0.82 & 1.00 & 1.00 & 0.35 & 0.42 & 0.50 & 0.74 & 0.80 & 0.81 \\ 
    & 50 & 0.69 & 1.00 & 1.00 & 0.10 & 0.13 & 0.21 & 0.72 & 0.74 & 0.79 \\ 
    & 100 & 0.53 & 0.90 & 1.00 & 0.08 & 0.15 & 0.15 & 0.56 & 0.74 & 0.79 \\ 
  Laplace & 20 & 0.86 & 1.00 & 1.00 & 0.93 & 1.00 & 1.00 & 0.81 & 0.83 & 0.81 \\ 
   & 50 & 0.72 & 1.00 & 1.00 & 0.55 & 1.00 & 1.00 & 0.89 & 0.91 & 0.91 \\ 
     & 100 & 0.52 & 0.89 & 1.00 & 0.19 & 0.89 & 1.00 & 0.82 & 0.95 & 0.94 \\   
  \hline
\end{tabular}
\caption{Mean of the sensitivity (Sen.: $\text{card}(\widehat{\Omega}\cap\Omega)/6$) for the feature selection obtained by the proposed method (Proposed method), the parametric method (VarSelLCM) and the sparse K-means (Sparcl) on 100 replicates for each scenario with theoretical misclassification rate of $5\%$, when the number of components is known.} \label{tab:simu1feature.sen}
\end{table}

\begin{table}[ht] 
\centering
\begin{tabular}{ccccccccccc}
  \hline 
& & \multicolumn{3}{c}{Proposed method}& \multicolumn{3}{c}{VarSelLCM} & \multicolumn{3}{c}{Sparcl} \\
& & \multicolumn{3}{c}{$n$}& \multicolumn{3}{c}{$n$} & \multicolumn{3}{c}{$n$} \\
Component & $J$ &  100 & 250 & 500& 100 & 250 & 500& 100 & 250 & 500 \\
  \hline
  Gaussian & 20 & 0.98 & 1.00 & 1.00 & 1.00 & 1.00 & 1.00 & 0.79 & 0.66 & 0.50 \\ 
     & 50 & 0.98 & 1.00 & 1.00 & 1.00 & 1.00 & 1.00 & 0.71 & 0.45 & 0.30 \\ 
   & 100 & 0.98 & 1.00 & 1.00 & 1.00 & 1.00 & 1.00 & 0.75 & 0.39 & 0.22 \\ 
    Student & 20 & 0.97 & 1.00 & 1.00 & 0.70 & 0.57 & 0.47 & 0.77 & 0.58 & 0.61 \\ 
    & 50 & 0.98 & 1.00 & 1.00 & 0.74 & 0.66 & 0.58 & 0.79 & 0.83 & 0.62 \\ 
    & 100 & 0.98 & 1.00 & 1.00 & 0.76 & 0.69 & 0.63 & 0.85 & 0.85 & 0.76 \\ 
   Laplace & 20  & 0.98 & 1.00 & 1.00 & 0.94 & 0.96 & 0.97 & 0.76 & 0.85 & 0.83 \\ 
   & 50 & 0.98 & 1.00 & 1.00 & 0.92 & 0.97 & 0.97 & 0.75 & 0.78 & 0.74 \\ 
  & 100  & 0.99 & 1.00 & 1.00 & 0.92 & 0.97 & 0.98 & 0.75 & 0.62 & 0.64 \\ 
  \hline
\end{tabular}
\caption{Mean of  the specificity (Spe.: $\text{card}(\widehat{\Omega}^c\cap\Omega^c)/(J-6)$) for the feature selection obtained by  the proposed method (Proposed method), the parametric method (VarSelLCM) and the sparse K-means (Sparcl)  on 100 replicates for each scenario with theoretical misclassification rate of $5\%$, when the number of components is known.} \label{tab:simu1feature.spe}
\end{table}

\paragraph{Full model selection}
We now compare both non-parametric and parametric approaches on their performances for full model selection. Table~\ref{tab:simu2K} presents the statistics on the number of components selected by both approaches. Again, when the distribution of the components is well-specified, the parametric approach obtains better results on  small samples because the proposed approach tends to underestimate the number of components. However, when the sample size increases, both methods perform similarly. When the distribution of the components is not Gaussian, the parametric method performs poorly and asymptotically overestimate the number of components with probability one. The proposed method is consistent for any number of variables, however, it tends to underestimate the number of components for small samples.

\begin{table}[ht]
\centering
\begin{tabular}{ccccccccccccccc}
  \hline  
Component & $J$ & \multicolumn{6}{c}{Proposed method}& \multicolumn{6}{c}{VarSelLCM} \\
& & \multicolumn{2}{c}{$n=100$}& \multicolumn{2}{c}{$n=250$}& \multicolumn{2}{c}{$n=500$}& \multicolumn{2}{c}{$n=100$}& \multicolumn{2}{c}{$n=250$}& \multicolumn{2}{c}{$n=500$}\\
& & Tr. & Ov. & Tr. & Ov.& Tr. & Ov.& Tr. & Ov.& Tr. & Ov.& Tr. & Ov. \\ 
  \hline
  Gaussian & 20  & 0.46 & 0.00 & 1.00 & 0.00 & 1.00 & 0.00 & 0.94 & 0.00 & 1.00 & 0.00 & 1.00 & 0.00 \\ 
     & 50 & 0.27 & 0.00 & 0.98 & 0.00 & 1.00 & 0.00 & 0.86 & 0.01 & 1.00 & 0.00 & 1.00 & 0.00 \\ 
     & 100 & 0.05 & 0.00 & 0.62 & 0.00 & 1.00 & 0.00 & 0.54 & 0.00 & 1.00 & 0.00 & 1.00 & 0.00 \\ 
    Student & 20 & 0.56 & 0.00 & 1.00 & 0.00 & 1.00 & 0.00 & 0.62 & 0.15 & 0.18 & 0.80 & 0.00 & 1.00 \\ 
     & 50 & 0.37 & 0.00 & 1.00 & 0.00 & 1.00 & 0.00 & 0.81 & 0.13 & 0.36 & 0.64 & 0.02 & 0.98 \\ 
    & 100 & 0.13 & 0.00 & 0.82 & 0.01 & 1.00 & 0.00 & 0.70 & 0.29 & 0.22 & 0.78 & 0.01 & 0.99 \\ 
      Laplace & 20 & 0.57 & 0.00 & 1.00 & 0.00 & 1.00 & 0.00 & 0.77 & 0.08 & 0.34 & 0.66 & 0.00 & 1.00 \\ 
     & 50 & 0.33 & 0.00 & 1.00 & 0.00 & 1.00 & 0.00 & 0.40 & 0.00 & 0.51 & 0.49 & 0.01 & 0.99 \\ 
   & 100 & 0.15 & 0.00 & 0.88 & 0.01 & 1.00 & 0.00 & 0.09 & 0.00 & 0.61 & 0.21 & 0.01 & 0.99 \\  
   \hline
\end{tabular}
\caption{Probability to select the true number of components (Tr.) and to overestimate it (Ov.) obtained by the proposed method (Proposed method) and the parametric method (VarSelLCM) on 100 replicates for each scenario  with theoretical misclassification rate of $5\%$, by performing a selection of  the variables.}\label{tab:simu2K}
\end{table}

Table~\ref{tab:simu2feature} presents the sensitivity and the specificity for feature selection obtained by both approaches when the number of components is also estimated. Again, results show the benefits of the proposed approach when the parametric assumptions are violated. In such a case, the parametric approach overestimates the number of components and, for heavy tail distributions (\emph{e.g.,} Student distribution), this approach tends to overestimate the subset of relevant variables. Moreover, for the small samples, the sensitivity is quite low explaining the tendency of overestimating the number of components. 
\begin{table}[ht]
\centering
\begin{tabular}{ccccccccccccccc}
  \hline  
Component & $J$ & \multicolumn{6}{c}{Proposed method}& \multicolumn{6}{c}{VarSelLCM} \\
& & \multicolumn{2}{c}{$n=100$}& \multicolumn{2}{c}{$n=250$}& \multicolumn{2}{c}{$n=500$}& \multicolumn{2}{c}{$n=100$}& \multicolumn{2}{c}{$n=250$}& \multicolumn{2}{c}{$n=500$}\\
& &  Sen. & Spe. & Sen. & Spe. & Sen. & Spe. & Sen. & Spe. & Sen. & Spe. & Sen. & Sep. \\ 
  \hline
 Gaussian & 20 & 0.97 & 0.82 & 1.00 & 1.00 & 1.00 & 1.00 & 1.00 & 0.98 & 1.00 & 1.00 & 1.00 & 1.00 \\ 
    & 50 & 0.95 & 0.74 & 1.00 & 0.99 & 1.00 & 1.00 & 1.00 & 0.96 & 1.00 & 1.00 & 1.00 & 1.00 \\ 
    & 100  & 0.94 & 0.59 & 1.00 & 0.94 & 1.00 & 1.00 & 0.99 & 0.87 & 1.00 & 1.00 & 1.00 & 1.00 \\ 
   Student & 20  & 0.96 & 0.86 & 1.00 & 1.00 & 1.00 & 1.00 & 0.67 & 0.37 & 0.66 & 0.72 & 0.56 & 0.97 \\ 
   & 50  & 0.96 & 0.79 & 1.00 & 1.00 & 1.00 & 1.00 & 0.73 & 0.09 & 0.68 & 0.15 & 0.59 & 0.24 \\ 
    & 100  & 0.95 & 0.70 & 1.00 & 0.98 & 1.00 & 1.00 & 0.76 & 0.09 & 0.72 & 0.13 & 0.67 & 0.13 \\ 
 Laplace & 20  & 0.97 & 0.87 & 1.00 & 1.00 & 1.00 & 1.00 & 0.92 & 0.90 & 0.96 & 1.00 & 0.91 & 1.00 \\ 
   & 50  & 0.97 & 0.80 & 1.00 & 1.00 & 1.00 & 1.00 & 0.86 & 0.44 & 0.97 & 1.00 & 0.96 & 1.00 \\ 
    & 100  & 0.95 & 0.72 & 1.00 & 0.99 & 1.00 & 1.00 & 0.84 & 0.13 & 0.96 & 0.83 & 0.98 & 1.00 \\
   \hline
\end{tabular}
\caption{Mean of the sensitivity (Sen.: $\text{card}(\widehat{\Omega}\cap\Omega)/6$) and the specificity (Spe.: $\text{card}(\widehat{\Omega}^c\cap\Omega^c)/(J-6)$) for the feature selection obtained by the proposed method (Proposed method) and the parametric method (VarSelLCM) on 100 replicates for each scenario with theoretical misclassification rate of $5\%$, when the number of components also is estimated.}\label{tab:simu2feature} 
\end{table}

\begin{figure}[ht!]
\centering \includegraphics[scale=0.36]{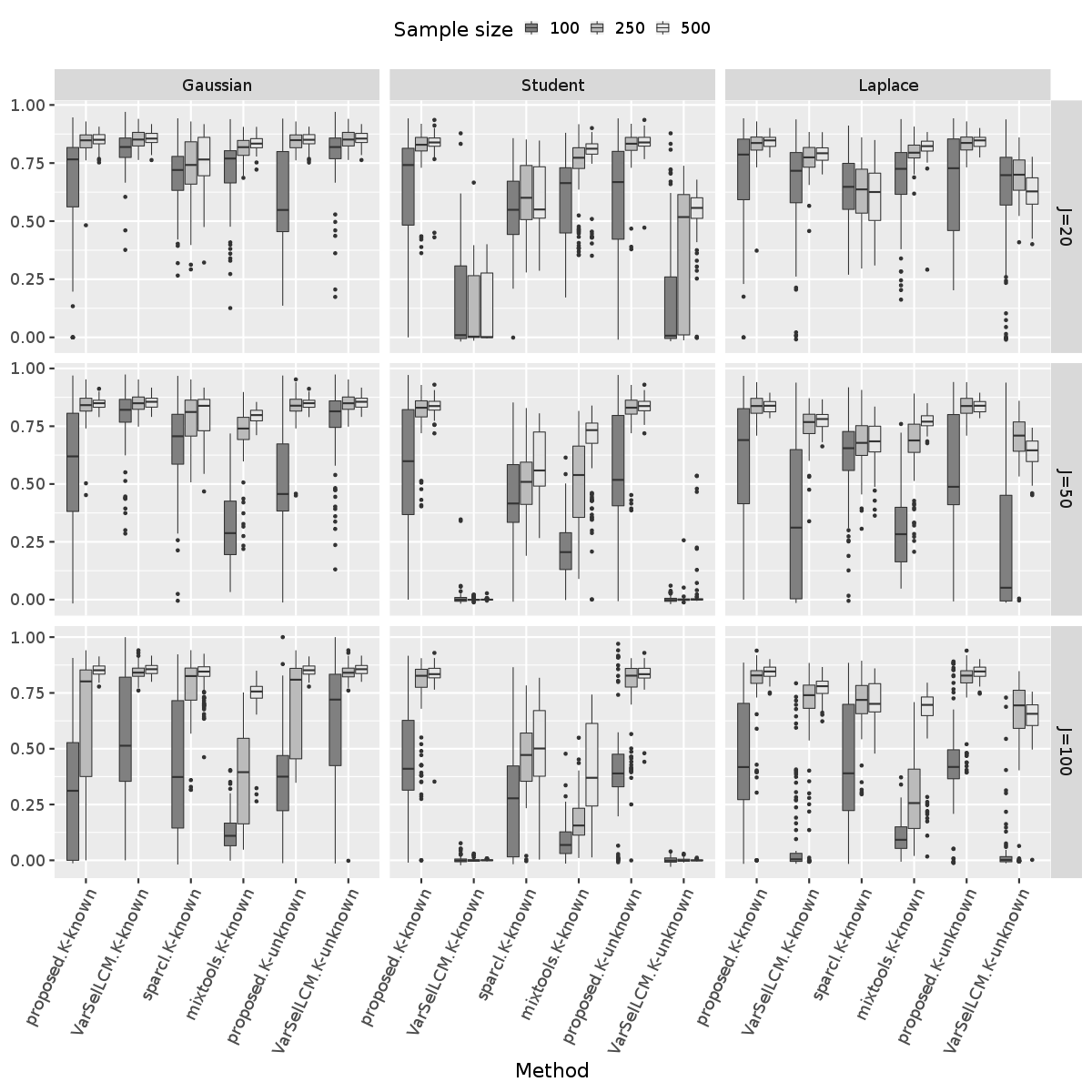}
\caption{Boxplot of the Adjusted Rand Index (ARI) obtained on the resulting partition when feature selection is performed with the true number of components by the proposed method (proposed.K-known) and by the parametric method (VarSelLCM.K-known), by the sparse K-means (Sparcl.K-known) and by the model considering all the variables as relevant components (mixtools.K-known) and when the full model selection (feature selection and estimation of the number of components) is achieved by the proposed approach (proposed.K-unknown) and the parametric approach (VarSelLCM.K-unknown). Data are generated with theoretical misclassification rate of $5\%$}\label{fig:resultsfeatureselection}
\end{figure}

\paragraph{Accuracy of the partition}
We are now interested in investigating the accuracy of the estimated partition. Thus, we compute the Adjusted Rand index \citep{HubertJoC1985} between the true partition and the estimators of the partition given by the non-parametric and the parametric methods when $K$ is known and then when it is estimated. Moreover, to illustrate the benefit of feature selection, we also estimate the partition by considering the full variables as relevant and the true number of components. Results are presented in Figure~\ref{fig:resultsfeatureselection}. Thus, when the parametric assumptions are satisfied, the parametric approach outperforms the proposed approach only on small samples (few observations with respect to the number of variables), whenever the number of components is known or not. However, when the parametric assumptions are violated, the proposed approach strongly outperforms the parametric approach. Note that, when the number of irrelevant variables increases, the approach considering all the variables for clustering performs poorly (see row 100), illustrating the benefit of feature selection for clustering.

\subsection{Benchmark data}
This section illustrates our procedure on four real data sets. The first data set illustrates the advantage of the procedure for selecting the number of components while the second data set sheds light on the importance of variable selection. The third data set shows that the procedure can be easily extended to the case of mixed-type data sets (a data set composed of continuous and categorical data). The last data set shows that the procedure can also be used to cluster data with non-ignorable missingness by considering the semi-parametric mixture models defined in \citet{DuRoy2020clustering}.

\paragraph{Swiss banknotes data}
We consider the Swiss banknotes data set \citep{Flury1988} containing six measurements (length of bill, width of left edge, width of right edge, bottom margin width, top margin width and length of diagonal) made on 100 genuine and 100 counterfeit old-Swiss 1000-franc bank notes. This data set is available in the R package mclust \citep{ScruccaR2016}. The status of the banknote (genuine or counterfeit) is also known. We perform the clustering of the bills based on the six morphological measurements and we evaluate the resulting partition with the status of the bills. Considering all the six morphological measurements as relevant for clustering, the proposed method detects two clusters which are strongly similar to the status of the bill (the ARI is 0.98 and only one genuine bill is assigned to the cluster grouping all the counterfeit bills). Clustering with Gaussian mixture models provides more components (mclust selects three components and VarSelLCM selects four components) and a partition related but different to the status of the bill (the ARI is 0.84 and 0.48 for mclust and VarSelLCM respectively). When a full model selection (feature selection and estimation of the number of components) is performed, the proposed method still selects two components and detects all the variables as relevant. Thus, a full model selection approach provides the same results as a method used for selecting the clusters by considering all the variables as relevant. Moreover, the Gaussian mixture models obtains less relevant results because VarSelLCM considers that all the measurements are relevant and thus obtains the same results as without performing feature selection.

\paragraph{Chemical properties of coffees}
We consider the data set collected by \citet{streuli1973heutige} that reports on the chemical composition of coffee samples collected from around the world. A total of 43 samples were collected from 29 countries, with beans from both Arabica and Robusta species, which is often considered as a pertinent partition. This data is available in the R package pgmm \citep{mcnicholas2015package}. We cluster the different coffee based on twelve  chemical constituents. A full Gaussian mixture clustering implemented in Mclust estimates three clusters and provided an ARI of 0.38. The same partition is obtained when the clustering is performed by VarSelLCM with a full model selection conducted according to the BIC (all the variables are detected as relevant for clustering). Again, similar results are obtained by the non-parametric mixture if the proposed method is used to select the number of components. However, if we perform a full model selection, only five of the twelve variables are detected as relevant for clustering and only two components are estimated. Moreover, this simpler model provides a perfect recovery of the species (ARI=1.00). This illustrates the importance of variable selection for clustering. Note that \citet{mcnicholas2008parsimonious} proposed a parsimonious Gaussian mixture model where parsimonious constraints are made on the covariance matrices of the Gaussian distributions, that also provides a perfect recovery of the partition.

\paragraph{Cleveland data set}
We consider the Cleveland dataset (available at \url{https://www.kaggle.com/ronitf/heart-disease-uci/version/1}). This data set is composed of $n=303$ subjects. Each subject is described by eight categorical variables having between two and six levels and five continuous variables. The ``goal" field refers to the presence of heart disease in the patient (no presence vs presence). Model \eqref{eq:modelavecselec1} can be easily extended to the case of mixed-type data (data set composed of continuous and categorical variables). Indeed, if variable $j$ is categorical then $\eta_{kj}$ is the probability mass function of a multinomial distribution. Thus, the discretization procedure used for model selection is applied only on the continuous variables while the number of levels for the categorical variables is fixed (\emph{i.e.,} it is not defined from the sample size). When the model is selected, the estimation of the extension of model \eqref{eq:modelavecselec1} can be easily achieved by maximizing the smoothed log-likelihood via an MM algorithm. The proposed approach detects the true number of clusters (\emph{i.e.,} two) while the approach implemented in VarSelLCM overestimates it since it selects six components. Moreover, by considering ten variables as relevant for clustering, our procedure returns a more relevant partition with respect to the occurrence of heart disease because it obtains an ARI equals to 0.37 while the procedure implemented in VarSelLCM obtains an ARI equal to 0.12.

\paragraph{Echocardiogram data set}
We consider the Echocardiogram Data Set \citep{salzberg1988exemplar} freely available in the R package MNARclust.  This data set is composed by $n=132$ subjects who suffered from heart attack at some point in the past. The task is generally to determine from the other variables whether or not the patient will survive at least one year. The data set is composed by 5 continuous variables: \emph{age at heart attack} (missing rate $4.5\%$), \emph{fractional shortening} (a measure of contracility around the heart, lower numbers are increasingly abnormal, missing rate $6.0\%$), \emph{epss} (E-point septal separation, another measure of contractility, larger numbers are increasingly abnormal, missing rate $11.4\%$), \emph{lvdd} (left ventricular end-diastolic dimension; this is a measure of the size of the heart at end-diastole; large hearts tend to be sick hearts, missing rate $8.3\%$) and \emph{wall-motion-score} (a measure of how the segments of the left ventricle are moving, missing rate $3.0\%$); one binary variable \emph{pericardial effusion} (pericardial effusion is fluid around the heart, 0=no fluid, 1=fluid, missing rate $0.7\%$). We also have one binary variable which can be used as a partition among the subjects:  \emph{still alive} (0=dead at end of survival period, 1=still alive). This binary variable is not used for clustering but permits to evaluate the accuracy of the estimated partition. 
Among the variables used for clustering there is 5.7\% of missing values and 19.1\% of the subjects have at least one missing value.  Moreover, the variable \emph{still alive} has only one missing value. \citet{DuRoy2020clustering} perform the cluster analysis of this data set by considering an extension of \eqref{eq:modelavecselec1} to the case of mixed-type data (data set composed of one binary and five continuous variables) with non-ignorable missingness process. Again, to use the proposed approach, we only discretize the continuous data. Moreover, to deal with the no-ignorable missingness, the absence of response is defined as a level of each variable. Performing a full model selection, the proposed approach detects three clusters. This results was suggested in \citet{DuRoy2020clustering} with empirical analysis of the evolution of the smoothed log-likelihood with respect to the number of clusters. Moreover, the proposed approach selects three relevant variables (fractional.shortening, epss and lvdd). Note that these three variables are detected as the most discriminative one for the missingness process and for their conditional distribution within component given the fact that the variable is observed (see Figures D.7 and D.8 in  \citet{DuRoy2020clustering}). Finally, note that the procedure returns a partition which is partially related to the indicator of surviving after the study since the ARI is 0.19.

\section{Conclusion}\label{sec:concl}
In this paper, we introduced a novel approach for a full model selection (number of components and subset of relevant variables for clustering) in multivariate finite mixture models. Thus, this approach is the first that allows variable selection for non-parametric mixture models. Under mild assumptions on the distributions of the observed variables, we showed that the model can be identified using a discretization of the data into bins and the penalized log-likelihood of the resulting distribution. Ranges  are given for the number of bins and the penalty term to ensure the consistency of the procedure. 
With a careful reading of our proof, a finite-sample size control of the probability of overestimating the model can be obtained.

 \bibliographystyle{chicago}
\bibliography{biblio}

\end{document}


\maketitle
\begin{abstract}
This supplement is organized as follows. Appendix A contains additional proofs, while Appendix B collects additional simulations.
\end{abstract}

\appendix

\section{Additional proofs} \label{sec:maths}

\subsection{Proof related to the underestimation} \label{proof:cvprobaunderestimation}
\subsubsection{Useful results}
\begin{lemma}[Hoeffding's Inequality]\label{lem:hoef}
Let $U_1,\ldots,U_n$ be independent real-valued random variables and $S_n$ be their sum $S_n= U_1+ \ldots + U_n$. Assume that there exist two real-valued sequences  $(a_k)_k$ and $(b_k)_k$ such that, for any $k$, $a_k < b_k$ and $\mathbb{P}\left(a_k < U_k < b_k \right)=1$. 
Then, for any positive $t$,
$$\mathbb{P}\left( S_n- \mathbb{E}[S_n] > t \right) \leq e^{- \frac{2 t^2}{\sum_{k=1}^n (b_k-a_k)^2}}.$$
\end{lemma}

\subsubsection{Proof of \eqref{eq:cvprobaunderestimation}}
Recall that we consider $M \in \mathcal{N}_1$. Note that $\Theta_{M,B, \varepsilon}\subset [\varepsilon,1]^{\text{dim(M,B)}}$ where $\text{dim(M,B)}=K + KB|\Omega| + B(J - |\Omega|)$. We denote by $\|\cdot \|_{\infty}$ the usual infinite norm on $\mathbb{R}^{\text{dim(M,B)}}$. Let the set
$$
\Gamma_{M,B,\rho}=\left\{\varepsilon + h\rho - \frac{\rho}{2}: h=1,\ldots,\left\lceil\frac{1}{\rho}\right\rceil \right\}^{\text{dim}(M,B)},
$$
where $\lceil\cdot\rceil$ is the ceiling operator. 
The set $[\varepsilon,1]^{\text{dim}(M,B)}$ can be covered by 
$\text{card}\left(\Gamma_{M,B,\rho}\right)=\left\lceil 1/\rho \right\rceil^{\text{dim}(M,B)}$ balls with diameter $\rho$ and whose  centers are the elements of $\Gamma_{M,B,\rho}$.  Therefore,
$\Theta_{M,B,\varepsilon}$ can be covered by 
$N_{M,B,\rho}:=\text{card}\left(\Gamma_{M,B,\rho}\right)$ balls with centers $\bvarpi_1,\ldots, \bvarpi_{N_{M,B,\rho}}$ and diameter $\rho$. 
Thus, defining $$m_{M,B,\rho}(\bx)=\sup_{\substack{ \theta, \tilde\theta\in\Theta_{M,B,\varepsilon}, \\ \|\theta - \tilde\theta\|_{\infty} < \rho,}} | \ln f_{M,B,\theta}(\bx) - \ln f_{M,B,\tilde\theta}(\bx)|,$$
we obtain that
\begin{equation}\label{eq:inequalityCVproba}
\sup_{\theta\in\Theta_{M,B,\varepsilon}} \frac{1}{n}\left(\ell_n(f_{M,B,\theta }) - \ell_n(g_0)\right)\leq \max_{s=1,\ldots,N_{M,B,\rho}} \frac{1}{n} \left(\ell_n(f_{M,B,\bvarpi_s} ) - \ell_n(g_0)\right) +\frac{1}{n} \sum_{i=1}^n m_{M,B,\rho}(\bx_i).
\end{equation}
We first deal with the last term of the right hand side of \eqref{eq:inequalityCVproba}, by showing  that
$\frac{1}{n} \sum_{i=1}^n m_{M,B,\rho}(\bx_i)=o_\mathbb{P}(1)$. Indeed, using that $|\ln a - \ln b|\leq |a-b|/(a \wedge b)$, we have
$$\left| \ln f_{M,B,\theta}(\bx_i) - \ln f_{M,B,\tilde\theta}(\bx_i)\right| \leq \frac{1}{K\varepsilon^{J+1}} \left| f_{M,B,\theta}(\bx_i) -  f_{M,B,\tilde\theta}(\bx_i)\right|.$$
Moreover, 
we have
\begin{multline}
\left| f_{M,B,\theta}(\bx_i) -  f_{M,B,\tilde\theta}(\bx_i) \right|
\leq  \sum_{k=1}^K \left|\pi_k - \tilde{\pi}_k\right| \prod_{j,b} \left(\frac{\alpha_{Bkjb}}{l_{Bjb}}\right)^{\sigma_{Bjb}(x_{ij})} +\\ \sum_{k=1}^K \tilde{\pi}_k \left|\prod_{j,b} \left(\frac{\alpha_{Bkjb}}{l_{Bjb}}\right)^{\sigma_{Bjb}(x_{ij})}-\prod_{j,b} \left(\frac{\tilde{\alpha}_{Bkjb}}{l_{Bjb}}\right)^{\sigma_{Bjb}(x_{ij})}\right|.
\end{multline}
On the one hand, for any $\bx \in\mathcal{X}$, any $\theta\in\Theta_{M,B,\varepsilon}$ and any $\tilde\theta\in\Theta_{M,B,\varepsilon}$ with $\|\theta - \tilde\theta\|_\infty < \rho$, we have  
$$ \sum_{k=1}^K \left|\pi_k - \tilde{\pi}_k\right| \prod_{j,b} \left(\frac{\alpha_{kjb}}{l_{Bjb}}\right)^{\sigma_{jb}(x_{ij})}  \leq  \frac{K \rho}{L_B^J},$$
where $L_B= \min_{j,b} l_{Bjb}$ is the size of the smallest bin.
On the other hand, using that, for any $|a_1| \leq 1$ and $|b_2| \leq 1$, $|a_1 a_2 - b_1 b_2| \leq |a_1- b_1| + |a_2-b_2|$, we obtain that:
$$
\sum_{k=1}^K \tilde{\pi}_k \left|\prod_{j,b} \left(\frac{\alpha_{kjb}}{l_{Bjb}}\right)^{\sigma_{jb}(x_{ij})}-\prod_{j,b} \left(\frac{\tilde{\alpha}_{kjb}}{l_{Bjb}}\right)^{\sigma_{jb}(x_{ij})}\right|    \leq \displaystyle  \sum_{k=1}^K \tilde{\pi}_k \frac{1}{L_B^J} \sum_j \left|\prod_{b} \alpha_{kjb}^{\sigma_{jb}(x_{ij})}-\prod_{b} \tilde{\alpha}_{kjb}^{\sigma_{jb}(x_{ij})}\right|.
$$
\color{blue}Noting that for any $x_{ij}$, there exists a single interval $b$ such that $\sigma_{jb}(x_{ij})=1$ (otherwise  $\sigma_{jb}(x_{ij})=0)$, we have $\left|\prod_{b} \alpha_{kjb}^{\sigma_{jb}(x_{ij})}-\prod_{b} \tilde{\alpha}_{kjb}^{\sigma_{jb}(x_{ij})}\right|\leq \rho$. \color{black}Thus,  for any $\bx \in\mathcal{X}$, any $\theta\in\Theta_{M,B,\varepsilon}$ and any $\tilde\theta\in\Theta_{M,B,\varepsilon}$ with $\|\theta - \tilde\theta\|_\infty < \rho$,
$$\sum_{k=1}^K \tilde{\pi}_k \left|\prod_{j,b} \left(\frac{\alpha_{kjb}}{l_{Bjb}}\right)^{\sigma_{jb}(x_{ij})}-\prod_{j,b} \left(\frac{\tilde{\alpha}_{kjb}}{l_{Bjb}}\right)^{ \sigma_{jb}(x_{ij})}\right| \leq \displaystyle \frac{J \rho}{L_B^J}.$$
This implies that, for any $\bx \in \mathcal{X}$, 
$$m_{M,B, \rho}(\bx) \leq  \frac{K+J}{K\varepsilon^{J+1}} \frac{\rho}{L_B^J}.$$
Thus, the radius $\rho$ will need to decrease at an appropriate rate as $n$ tends to infinity in order to satisfy the announced convergence. Note that by definition $\varepsilon^{-1}=O(n^{\alpha + 1})$. Moreover, using Assumption~\ref{ass:intervals}, one can find some positive constant $\gamma$ such that $(L_B^J)^{-1}=O(n^{J\gamma})$. Thus, we set $\rho = o(n^{-((J+1)(\alpha + 1) + J\gamma)})$ to obtain
\begin{equation}\label{eq:under1}
\frac{1}{n} \sum_{i=1}^n m_{M,B,\rho}(\bx_i)  = o_\mathbb{P}(1).
\end{equation}
Now we investigate the asymptotic behavior of the first term in the right-hand side of  \eqref{eq:inequalityCVproba} considering the decrease rate of $\rho$. We split it into three terms as follows
\begin{equation}\label{eq:under2}
\max_{s} \frac{1}{n} \left(\ell_n(f_{M,B,\bvarpi_s} ) - \ell_n(g_0)\right) \leq \max_{s} \Delta_{n,s} + \max_{s} \mathbb{E}_{g_0}[\ln f_{M,B,\bvarpi_s}-\ln g_0]  + \left( \mathbb{E}_{g_0}[\ln g_0] - \ell_n(g_0) \right),
 \end{equation}
where $\Delta_{n,s}=\frac{1}{n} \ell_n(f_{M,B,\bvarpi_s} ) - \mathbb{E}_{g_0}[\ln f_{M,B,\bvarpi_s}]$.  
First of all, the law of  large numbers permits to deal directly with the third term, as it implies that 
$$ \frac{1}{n} \ell_n(g_0) - \mathbb{E}_{g_0}[\ln g_0] = o_\mathbb{P}(1).$$
For the first term in \eqref{eq:under2}, we have, for any $\eta>0$,
$$
\mathbb{P}\left(\max_{s=1,\ldots,N_{M,B,\rho}} \Delta_{n,s}>\eta\right) \leq \sum_{s=1}^{N_{M,B,\rho}}\mathbb{P}\left(\Delta_{n,s} >\eta \right).
$$
For any $\bx$ and any $s$, we have the following bounds on $\ln f_{M,B,\bvarpi_s}$, 
$$ \ln K + (J+1) \ln \varepsilon \leq \ln f_{M,B,\bvarpi_s}(\bx)\leq \ln K - J \ln L_B.$$
Thus, we can apply Lemma~\ref{lem:hoef} to obtain an upper bound of $\mathbb{P}\left(\Delta_{n,s} >\eta \right)$. Thus, by noting that $N_{M,B,\rho}=\exp (\text{dim}(M,B) \ln\left\lceil 1/\rho\right\rceil )$, we have  for any $\eta>0$
$$
\displaystyle \mathbb{P}\left(\max_{s=1,\ldots,N_{M,B,\rho}} \Delta_{n,s}>\eta\right) 
\leq  \exp \left( -\frac{2 \eta^2 \; n}{(J \ln L_B + (J+1) \ln \varepsilon)^2 }  +\text{dim}(M,B) \ln\left\lceil\frac{1}{\rho}\right\rceil \right).
$$
Using Assumption~\ref{ass:intervals}.\ref{ass:intervals2} and the definition of $\varepsilon$, we have $\frac{2 \eta^2 \; n}{(J \ln L_B + (J+1) \ln \varepsilon)^2 } = O(n/\ln^2 n)$. Moreover,  the rate of decreasing of $\rho$ is such that $\ln\left\lceil 1 / \rho\right\rceil = O(\ln n)$. Finally, noting that $\text{dim}(M,B)=O(B)$, 
Assumption~\ref{ass:intervals}.\ref{ass:vitesseB} implies that
\begin{equation}\label{eq:under3}
\max_{s=1,\ldots,N_{M,B,\rho}} \Delta_{n,s}  = o_\mathbb{P}(1).
\end{equation}

The second term in \eqref{eq:under2} is the one that will not be asymptotically negligeable in probability. As the parameter space $\Theta_{M,B,\varepsilon}$ is compact, there exists some $\theta_B^\star \in \Theta_{M,B,\varepsilon}$ which maximizes the value of $\mathbb{E}_{g_0}[\ln f_{M,B,\theta}]$ over $\theta$. Straigthforwardly, we have that, for any $s$, 
 $$\mathbb{E}_{g_0}[\ln f_{M,B,\bvarpi_s}-\ln g_0] \leq  \mathbb{E}_{g_0}[\ln f_{M,B,\theta_B^\star}-\ln g_0].$$
Besides, $f_{M,B,\theta_B^\star}$ can be approximated by a function $\tilde g_{M,B}\in\mathcal{G}_M$ in such a way that $\| f_{M,B,\theta_B^\star} - \tilde g_{M,B}\|_{\infty}=O(n^{-\gamma J})$.
 Noting that $\text{KL}(g_0,\tilde g_{M,B})\geq \delta_M$, we deduce that $$ \mathbb{E}_{g_0}[\ln f_{M,B,\theta_B^\star}-\ln g_0] \leq -\delta_M + o_{\mathbb{P}}(1).$$
Thus, combining these results and \eqref{eq:under1}-\eqref{eq:under3} to control the terms in \eqref{eq:inequalityCVproba}, we obtain the following upper-bound
$$\sup_{\theta\in\Theta_{M,B,\varepsilon}} \frac{1}{n}\left(\ell_n(f_{M,B,\theta }) - \ell_n(g_0)\right)\leq -\delta_M + o_{\mathbb{P}}(1).$$

\subsection{Proof of \eqref{eq:cvover}} \label{sec:cvover}
To show that \eqref{eq:cvover} holds true, we use the locally conic parametrization introduced by \citet{DacunhaESAIM1997,DacunhaAOAS1999} that permits, in particular, to state the convergence in distribution of the log-likelihood ratio (see  Section 2.2 in \citet{KeribinSan2000}). The locally conic parametrization is introduced in Section~\ref{sec:loc}, then Section~\ref{sec:assDacunha} presents the assumptions required to obtain the convergence in distribution. Finally, Section~\ref{sec:bracketing} presents usefull tools about bracketing entropy and Section~\ref{sec:assok} shows that the assumptions we made on the data distribution imply that the assumptions detailled in Section~\ref{sec:assDacunha} are satisfied.

\subsubsection{Locally conic parametrization}\label{sec:loc}
The principle of the local parametrization is to circumvent the issue of non-identifiability of the parameters of model $M$ that contains the true model $M_0$ (\emph{i.e.,} $M_0\in M$). The idea is to consider a new parametrization for the mixture model \eqref{eq:modeldiscret} defined by
$$
f_{M,B,\theta}=\sum_{k=1}^K \pi_k f_{k}(x_{i};\alpha_k)
\; \text{ with } 
\; 
f_{k}(x_{i};\alpha_k)=\prod_{j=1}^J \prod_{b=1}^B \left(\frac{\alpha_{kjb}}{l_{jb}}\right)^{\sigma_{jb}(x_{ij})},
$$
where $M$ imposes that $\alpha_{1j}=\ldots=\alpha_{Kj}$ for any $j\in\bar{\Omega}$.
The new parametrization is such that a scalar positive parameter $\mu$ is identifiable and a second parameters $\beta$ contains all the non-identifiable parameters. Thus, we denote by $\phi_{M,B,\mu,\beta}$ the pdf of the mixture model obtained with the locally conic parametrization and we have
$$
\phi_{M,B,\mu,\beta}=f_{M_0,B,\theta_0} \Leftrightarrow \mu=0.
$$
Thus, $\mu$ measures the distance to the real model and $\beta$ represents the directions for which the sub-models approach the true model. For any model $M\in\mathcal{N}_3$, $\mathcal{B}_{K,B}$ is the set of parameters 
$$\beta=(\lambda_1,\ldots,\lambda_{K-K_0},\alpha_{11}^\top,\ldots,\alpha_{(K-K_0) J}^\top,\delta_{11}^\top,\ldots,\delta_{K_0 J}^\top,\rho_1,\ldots,\rho_{K_0})^\top,$$
where, for $\ell=1,\ldots,K-K_0$, $k=1,\ldots,K_0$ and $j=1,\ldots J$, $\lambda_\ell\geq 0$, $\alpha_{\ell j}=\mathcal{S}_B$, $\delta_{kj} \in\mathbb{R}^{B}$ and  $\rho_k\in\mathbb{R}$. To ensure that the new parametrization is locally conic, restrictions have to be done on the parameter space $\mathcal{B}_{K,B}$ (see Section~3 in \citep{DacunhaAOAS1999}). Moreover, supplementary constraints have to be considered to respect the role of the variables defined by $\Omega$ (relevant or irrelevant for clustering). Thus, the resulting space of $\beta$ is denoted by $\tilde{\mathcal{B}}_{K,B}$ and the space of $(\mu, \beta^\top)^\top$ is denoted by $\mathcal{T}_B=[0,R]\times\tilde{\mathcal{B}}_{K,B}$. The locally conic parametrization can be understood as a perturbation of the true distribution $f_{0,B}$ defined as perturbation of the parameters $\alpha_{0,k}$ and the weights $\pi_{0,k}$ for the  $K_0$ mixture $f_{0,B}$ and an addition of a $K-K_0$-mixture whose weights tend to zero. Thus, the pdf obtained with the locally conic parametrization is
$$
\phi_{M,B,\mu,\beta}=\sum_{\ell=1}^{K-K_0} \lambda_\ell \frac{\mu}{N(\beta)}\varphi_{\alpha_{\ell}} + \sum_{k=1}^{K_0} \left(\pi_{0,k} + \rho_k \frac{\mu}{N(\beta)}\right) \varphi_{\alpha_{0,k} + \frac{\mu}{N(\beta)} \delta_k},
$$
with
$$
\varphi_{\upsilon}(\bx_i) = \prod_{j=1}^J \prod_{b=1}^B \left( \frac{\upsilon_{jb}}{l_{jb}} \right)^{\sigma_{jb}(x_{ij})},
$$
and 
$$
N(\beta) = \left\| \sum_{k=1}^{K_0} \pi_{0,k} \sum_{j=1}^{J}\sum_{b=1}^B \delta_{kjb} \frac{D_{jb}^1 \varphi_{\alpha_{0,k}}}{f_{0,B}} + \sum_{\ell=1}^{K-K_0} \lambda_{\ell} \frac{\varphi_{\alpha_{\ell}} }{f_{0,B}} + \sum_{k=1}^{K_0} \rho_k \frac{ \varphi_{\alpha_{0,k}}}{f_{0,B}} \right\|_{L^2_{f_{0,B}}},
$$
where $D_{jb}^1 \varphi_{\alpha_{0,k}}$ is the derivative of $ \varphi_{\alpha_{k}}$ with respect to the $b$-th element of vector $\alpha_{kj}$ and evaluated at point $\alpha_k=\alpha_{0,k}$. If the assumptions stated in Section~\ref{sec:assDacunha} hold true then equation \eqref{eq:cvover} is satisfied, and additionnally, for a fixed value $B^\star$, the likelihood ratio converges in distribution as follows
\begin{equation} 
T_{n,M,B^\star} - \ell_n f_{0,B^\star} \xrightarrow{d} \sup\left\{\sup_{d\in \mathcal{D}_{B^\star}}\frac{1}{2} \xi_d^2 \mathds{1}_{\xi_d\geq 0}, \sup_{d_1\in \mathcal{D}_{1,B^\star}; d_2\in \mathcal{D}_{2,B^\star}}\frac{1}{2} \left(\xi_{d_1}^2  + \xi_{d_2}^2 \mathds{1}_{\xi_{d_2}\geq 0}\right)\right\},
\end{equation}
where $\xi_d$, $\xi_{d_1}$ and $\xi_{d_2}$ are the Gaussian processes indexed by $\mathcal{D}_{B}$, $\mathcal{D}_{1,B}$ and $\mathcal{D}_{2,B}$ with covariance the usual Hilbertian product in $L^2_{f_{0,B}}$, $\mathcal{D}_{B}$ has a compact closure $\overline{\mathcal{D}}_B$ that is the subset of the unite sphere of $L^2_{f_{0,B}}$ of functions of form
$$
\frac{1}{N(\beta)} \left(\sum_{k=1}^{K_0} \pi_{0,k} \pi_{0,k} \sum_{j=1}^{J}\sum_{b=1}^B \delta_{ks} \frac{D_{jb}^1 \varphi_{\alpha_{0,k}}}{f_{0,B}} + \sum_{\ell=1}^{K-K_0} \lambda_{\ell} \frac{\varphi_{\alpha_{\ell}} }{f_{0,B}} + \sum_{k=1}^{K_0} \rho_k \frac{ \varphi_{\alpha_{0,k}}}{f_{0,B}} \right),
$$
with $\beta \in \tilde{\mathcal{B}}_{K,B}$ and where $\mathcal{D}_{1,B}$ and $\mathcal{D}_{2,B}$ are two orthogonal subsets of $\mathcal{D}_B$, fully described in \citet{KeribinSan2000}.

\subsubsection{Assumptions of \citet{DacunhaAOAS1999}}\label{sec:assDacunha}

\textbf{(P0)} There exists a function $\omega$ in $L^1_{f_{0,B}}$ such that $\forall \alpha_k \in \mathcal{S}_B^J$ $|\ln \varphi_{\alpha_k}|\leq \omega$ almost surely. Moreover, $ \varphi_{\alpha_k}$ possesses partial derivatives up to order 5. For all $q\leq 5$ and all couples $(j_1,b_1),\ldots, (j_q,b_q)$
$$
\frac{D^q_{j_1 b_1,\ldots,j_q b_q}  \varphi_{\alpha_{0,k}}}{f_{0,B}} \in L^3_{f_{0,B}}.
$$
where $D^q_{j_1 b_1,\ldots,j_q b_q}  \varphi_{\alpha_k}$ is the $q$-th derivative of $\varphi_{\alpha_k}$ with respect to $\alpha_{k j_1 b_1}, \ldots, \alpha_{k j_q b_q}$.
Moreover, there exists a function $m_5$ and a positive $\varepsilon$ such that
$$
\sup_{\|\alpha_k - \alpha_{0,k}\|\leq \varepsilon}
\left| \frac{D^5_{j_1 b_1,\ldots,j_5 b_5}  \varphi_{\alpha_k}}{f_{0,B}}  \right| \leq m_5
\;
\text{ and }
\;
\mathbb{E}_{f_{0,B}}[m_5^3] < \infty.
$$

\noindent \textbf{(P1)} For any integer $p_1$ and $p_2$ such that $p_1+p_2\leq K-K_0$, for any set of distinct points $\alpha_{\ell}$, $\ell=1,\ldots,p_1$, distinct from any $\alpha_{0,\ell}$, any permutation $\mathfrak{s}$ of $[1,\ldots,K_0]$, the set of functions
$$
\left(
\left(\frac{ \varphi_{\alpha_{\ell}}}{f_{0,B}}\right)_{\ell=1,\ldots,p_1},
\left(\frac{ \varphi_{\alpha_{0,{\ell}}}}{f_{0,B}}\right)_{\ell=1,\ldots,K_0},
\left(\frac{ D_{jb}^1\varphi_{\alpha_{0,{\ell}}}}{f_{0,B}}\right)_{\underset{j=1,\ldots,J; b=1,\ldots,B}{\ell=1,\ldots,K_0;} },
\left(\frac{ D_{j_1 b_1, j_2 b_2}^2\varphi_{\alpha_{0,{\ell}}}}{f_{0,B}}\right)_{\underset{\underset{b_1,b_2 = 1, \ldots,B}{j_1, j_2 = 1,\ldots,J;}}{\ell=\mathfrak{s}(1),\ldots,\mathfrak{s}(p_2);}}
\right)
$$
are linearly independent in $L^2_{f_{0,B}}$.

\subsubsection{Bracketing and entropy number}\label{sec:bracketing}
A relatively simple way to measure the size of any set is to use covering numbers. Let $(\mathcal{F},\|\cdot\|)$ be an arbitrary semi-metric space. A $\varepsilon$-cover of the set $\mathcal{F}$ with respect to the semi-metric $\|\cdot\|$ is the set of $(f_1,\ldots,f_N)$ (note that these elements need not belongs to $\mathcal{F}$ themselves) such that for any $f\in\mathcal{F}$, there exists some $v \in \{1,\ldots, N\}$ with $\| f - f_v\|\leq \varepsilon$. The $\varepsilon$-covering number $\mathcal{N}(\varepsilon,\mathcal{F},\|\cdot\|)$ is the minimal number of balls $B(g;\varepsilon)=\{f\in\mathcal{F}: \|f-g\|\leq \varepsilon\}$ of radius $\varepsilon$ needed to cover $\mathcal{F}$.

An alternative way to measure the size of any set is to use bracketing number. Given two functions $l$ and $u$, the bracket $[l,u]$ is the set of all the function $f\in\mathcal{F}$ with $l(x)\leq f(x) \leq u(x)$, for all $x\in\mathcal{X}$. An $\varepsilon$-bracket is a bracket $[l,u]$ with $\|l-u\|< \varepsilon$. The bracketing number $\mathcal{N}_{[]}(\varepsilon,\mathcal{F},\|\cdot\|)$ is the minimum number of $\varepsilon$-brackets needed to cover $\mathcal{F}$. Note that $l$ and $u$ need not belong to $\mathcal{F}$ themselves, but are assumed to have finite norms. $\mathcal{N}_{[]}(\varepsilon,\mathcal{F},\|\cdot\|)$ being equal to 1 for a sufficiently large $\varepsilon$, we only have to check the order of $\mathcal{N}_{[]}(\varepsilon,\mathcal{F},\|\cdot\|)$ when $\varepsilon$ tends to 0.

Note that there is the following relation between covering and bracketing numbers, for any $\varepsilon>0$:
\begin{equation}\label{eq:relationcoveringbracketing}
\mathcal{N}(\varepsilon, \mathcal{F},\|\cdot\|) \leq 
\mathcal{N}_{[]}(2\varepsilon, \mathcal{F},\|\cdot\|).
\end{equation}

The following lemma shows that, under mild conditions, there is a relation between the number of pavements needed to cover $\mathcal{B}_B$ and the bracketing number needed to cover $\mathcal{D}_B$. Its proof is given in Section 4.1 in \citet{KeribinSan2000}. To state the lemma, we need to introduce the following condition

\textbf{(P2)} All the products of functions in $\left\lbrace \left(\frac{ \varphi_{\alpha_{k}}}{f_{0,B}}\right),
\left(\frac{ D_{jb}^1\varphi_{\alpha_{k}}}{f_{0,B}}\right),
\left(\frac{ D_{j_1 b_1, j_2 b_2}^2\varphi_{\alpha_{k}}}{f_{0,B}}\right)\right\rbrace$ are square integrable for $f_{0,B}$.

\begin{lemma}[Number of pavements and bracketing number] \label{lemma:pavandbracke}
Assume that (P0), (P1) (P2) hold true. If we divide each set of parameters of $\tilde\beta_B$ in intervals of length $\varepsilon$, we obtain a partition of $\mathcal{D}_B$ of cardinal of order $O(\varepsilon^{- \tau})$ and then
$$
\mathcal{N}_{[]}(\varepsilon, \mathcal{D}_B,\|\cdot\|) = O(\varepsilon^{-\tau}).
$$
\end{lemma}

\subsubsection{Validation of the Assumptions of \citep{DacunhaAOAS1999}} \label{sec:assok}

First of all,  for any $q \leq 5$ and any subset $s_q = \{ (j_1 , b_1), \ldots, (j_q , b_q) \}$ of 
$\left(\{1, \ldots,J\} \times \{1, \ldots, B\}\right)^q$, we have the following derivative
$$D^q_{j_1 b_1,\ldots,j_q b_q}  \varphi_{\alpha_{k}} (\bx) = \prod_{(j,b) \in s_q} \frac{\sigma_{jb}(x_j)}{l_{jb}} \prod_{(j,b) \notin s_q} \left(\frac{\alpha_{kjb}}{l_{jb}} \right)^{\sigma_{jb}(x_j)},$$
if all elements in $s_q$ are distinct, and $D^q_{j_1 b_1,\ldots,j_q b_q}  \varphi_{\alpha_{k}}$ is equal to zero otherwise.
Therefore, (P1) is satisfied and it suffices to derive bounds for $f_{0,B}$ and for $\varphi_{\alpha_k}$ to obtain (P0) and (P2).
For any $\alpha_k \in S_B^J$ and for any $\bx \in \mathcal{X}$, $$0 \leq \varphi_{\alpha_k} (\bx) \leq L_B^{-J},$$  $$0 \leq D^q_{j_1 b_1,\ldots,j_q b_q}  \varphi_{\alpha_{k}} (\bx) \leq L_B^{-J}$$ and $$(\min_{l=1, \ldots, K} \pi_l) \varphi_{\alpha_{0,k}}(\bx) \leq f_{0,B}(\bx) \leq K L_B^{-J}.$$
 Additionnally, at point $\alpha_{0,k}$, by construction, 
$\alpha_{0,kjb} = \int_{I_{Bjb}} \eta_{kj}(x_j) \mathrm{d}x_j.$
Let $$\varepsilon_B = \min_{\underset{k=1,\ldots, K}{\underset{b=1, \ldots, B}{j=1, \ldots, J}}} \int_{I_{Bjb}} \eta_{kj}(x_j) \mathrm{d}x_j, $$ as the densities $\eta_{kj}$ cannot be null over any whole bin $I_{Bjb}$, we obtain that $\varepsilon_B >0$.
Thus, we deduce that, for any $k, j, b$, $\alpha_{0,kjb} \geq \varepsilon_B >0$, which implies that
$$f_{0,B} (\bx) \geq (\min_l \pi_l) \varepsilon_B^J \prod_{j,b} \left(\frac{1}{l_{jb}}\right)^{\sigma_{jb}(x_j)},$$ and that 
$$ 0 \leq \frac{D^q_{j_1 b_1,\ldots,j_q b_q}  \varphi_{\alpha_{k}}}{f_{0,B}}(\bx) \leq \left((\min_l \pi_l) \varepsilon_B^J\right)^{-1}.$$

\subsection{Proof of \eqref{eq:asymptV}} \label{sec:asymptV}
This section explains the control of the deviation of the supremum of the Gaussian process indexed by function on space $D_{B,s}$.
Hence we start by reminding the Dudley's inequality that control the expectation of the supremum of a Gaussian process and the Borell-TIS inequality that control the deviation for its mean of the supremum of a Gaussian process.

\begin{lemma}[Dudley's inequality]
Let $(\xi_t)$ be a Gaussian process indexed by a semi-metric space $(\mathcal{F},\|\cdot\|)$. For $\varepsilon>0$, denote by $ \mathcal{N}(\varepsilon, \mathcal{F},\|\cdot\|)$ the covering number, then there exists some universal constant $\kappa>0$ such that
$$
\mathbb{E}\left[\sup_{t\in \mathcal{F}} \xi_t\right] \leq \kappa \int_0^{\text{diam}(\mathcal{F},\|\cdot\|)/2} \sqrt{\ln \mathcal{N}(\varepsilon,\mathcal{F},\|\cdot\|)} d\varepsilon.
$$
\end{lemma}

\begin{lemma}[Borell-TIS inequality]
Let $(\xi_t)$ be a centered Gaussian process indexed by a semi-metric space $(\mathcal{F},\|\cdot\|)$. Let $\varsigma_\xi^2=\sup_{t\in \mathcal{F}} \text{Var}(\xi_t)$ be the supremum of the  variance of the process. If $(\xi_t)$ is a separable process then for any $\varepsilon>0$
$$
\mathbb{P}\left( \left| \sup_{t\in \mathcal{F}} |\xi_t| - \mathbb{E}\left[\sup_{t\in \mathcal{F}} |\xi_t| \right] \right| \geq \varepsilon \right) \leq 2 \exp\left( - \frac{\varepsilon^2}{2 \varsigma_\xi^2} \right).
$$
\end{lemma}

To prove \eqref{eq:asymptV}, we want to control, for any $\varepsilon>0$
\begin{equation} \label{eq:boudV1}
\mathbb{P}\left(\sup_{d\in \mathcal{D}_{B,s}} \xi_d \geq \sqrt{\varepsilon a_{n,M_0,B}}\right)
\leq \mathbb{P}\left(\left| \sup_{d\in \mathcal{D}_{B,s}} |\xi_d|  - \mathbb{E}\left[\sup_{d\in \mathcal{D}_{B,s}} |\xi_d| \right] \right| \geq \sqrt{\varepsilon a_{n,M_0,B}} - \mathbb{E}\left[\sup_{d\in \mathcal{D}_{B,s}} |\xi_d| \right]\right).
\end{equation}
In Section~\ref{sec:assok}, we shown that assumptions (P0), (P1) and (P2) hold true under the assumption of Theorem~\ref{thm:cvproba}. Therefore, applying Lemma~\ref{lemma:pavandbracke} and the relation between covering and bracketing numbers stated by \eqref{eq:relationcoveringbracketing}, and noting that we can divide $\tilde{\mathcal{B}}_B$ with a number of pavements of order $O(\varepsilon^{- B JK_{\max}})$, we have $\mathcal{N}(\varepsilon, \mathcal{D}_{B,s},\|\cdot\|_{L^2f_{0,B}}) = O((\varepsilon/2)^{-BJK_{\max}})$. Noting that $\text{diam}(\mathcal{D}_{B,s},\|\cdot\|_{L^2f_{0,B}})=1$, Dudley's inequality implies that there exists a positive constant $\tilde\kappa$ such that
\begin{equation} \label{eq:boudexpect}
\mathbb{E}\left[\sup_{d\in \mathcal{D}_{B,s}} |\xi_d| \right] \leq \tilde\kappa\sqrt{BJK_{\max}}.
\end{equation}
The space $\mathcal{D}_{B,s}$ is separable, thus,  replacing the expectation of the supermum in \eqref{eq:boudV1} by its upperbound defined by \eqref{eq:boudexpect} then applying the Borell-TIS inequality leads to
$$
\mathbb{P}\left(\sup_{d\in \mathcal{D}_{B,s}} \xi_d \geq \sqrt{\varepsilon a_{n,M_0,B}}\right)
\leq
2 \exp\left( - \frac{\sqrt{\varepsilon a_{n,M_0,B}} - \tilde\kappa\sqrt{(K_{\max}-1) + K_{\max} B J}}{2 \varsigma_\xi^2} \right).
$$
Assumption~\ref{ass:penalty} made on the penalty ensures that $\lim_{n\to \infty}  a_{n,M_0,B} / \left(K_{\max}-1) + K_{\max} B J\right) =0$ and thus
$$
\lim_{n\to\infty} \mathbb{P}\left(\sup_{d\in \mathcal{D}_{B,s}} \xi_d \geq \sqrt{\varepsilon a_{n,M_0,B}}\right) = 0,
$$
which concludes the proof of \eqref{eq:asymptV}.

\subsection{Proof of \eqref{eq:asymptconcentre}} \label{sec:asymptconcentre}
First, we presents the results on the Gaussian approximation of suprema of empirical processes stated in \citet{ChernozhukovAOS2014}. Then, we use those results to prove \eqref{eq:asymptconcentre}.

\subsubsection{Gaussian approximation of suprema of empirical processes}
All the elements presented in this section arise from \citet{ChernozhukovAOS2014}.
Let $\mathcal{U}$ be a class of measurable functions $\mathcal{X}\to\mathbb{R}$ that is $P$-centered (\emph{i.e.,} $Pu=0$ with $Pu=\int u dP$). Denote by $U$ a measurable envelope of $\mathcal{U}$, that is, $U$ is a non-negative measurable function $\mathcal{X}\to\mathbb{R}$ such that $U(x)\geq \sup_{u\in\mathcal{U}} |u(x)|$, $\forall x\in\mathcal{X}$. We make the following assumptions:
\begin{itemize}
\item[(C1)] The class $\mathcal{U}$ is point-wise measurable, that is, it contains a countable subset $\mathcal{V}$ such that for every $u\in\mathcal{U}$ there exists a sequence $v_m\in\mathcal{V}$ with $v_m(x)\to u(x)$ for every $x\in\mathcal{X}$.
\item[(C2)] For some $q\geq 4$, $U\in\mathcal{L}^q(P)$.
\item[(C3)] the class $\mathcal{U}$ satisfies for all $0<\varepsilon\leq 1$,
$$
\mathcal{N}_{[]}(\varepsilon \|U\|_{L^2P}, \mathcal{U}, \|\cdot\|_{L^2P}) = O\left( \left[\frac{A}{\varepsilon}\right]^{\nu} \right),
$$
where $A \geq e$ and $\nu>1$.
\end{itemize}
Note that the assumption (C1) is the Assumption (A1) in \citet{ChernozhukovAOS2014} and that Assumptions~(C2) and (C3) implies that the Assumptions (A2) and (A3)  in \citet{ChernozhukovAOS2014} hold true. Thus, Theorem 5.1 in \citet{ChernozhukovAOS2014} holds true under (C1)-(C3) but implies to define bounds on triple terms (see (5) in \citet{ChernozhukovAOS2014}). The authors control these terms using the uniform entropy integral (see Lemma 2.2 in \citet{ChernozhukovAOS2014}) but the same control can be obtain with the integral of the non-uniform entropy of bracketing. Moreover, a simple bound can be obtained to VC type classes of function. Despite that Assumption (C3) does not implies that $\mathcal{F}$ is a VC type class (because control of the entropy with bracketing is not uniform on $P$), Corollary 2.2  in \citet{ChernozhukovAOS2014} still holds under $(C3)$ (the proof is similar). Hence, we obtain the following lemma that control the Gaussian approximation of suprema of empirical processes.

\begin{lemma}[Gaussian approximation to suprema of empirical processes]\label{lem:empriproce}
Suppose that assumptions (C1)-(C3) are satisfied. In addition, suppose that for some $b\geq \sigma \geq 0$ and $q\in[4,\infty)$, we have $\sup_{u\in\mathcal{U}} P|u|^k\leq\sigma^2 b^{k-2}$ for $k=2,3$ and $\|U\|_{P,q}\leq b$. Let $Z=\sup_{u\in\mathcal{U}} \mathcal{G}_n u$ and $\tilde Z $ which follows in distribution $\sup_{u\in\mathcal{U}} G_p u$, then for every $\gamma \in (0,1)$
$$
\mathbb{P}\left[
| Z - \tilde Z | >
\frac{bK_n}{\gamma^{1/2}n^{1/2-1/q}} +
\frac{(b\sigma)^{1/2}K_n^{3/4}}{\gamma^{1/2}n^{1/4}} + 
\frac{(b\sigma^2 K_n^2)^{1/3}}{\gamma^{1/3}n^{1/6}}
\right]
\leq C\left(\gamma + \frac{\ln n}{n}\right),
$$
where $K_n=c\nu(\ln n \vee \ln (Ab/\sigma))$ and $c$ and $C$ are strictly positive constants that depend only on $q$ ("$1/q$" is interpreted as 0 when $q=\infty$).
\end{lemma}

\subsubsection{Application to the control of the log-likelihood ratio}
We use the notation of the previous lemma. In our framework, $Z=\sup_{d\in \mathcal{D}_{B,s}} \mathcal{G}_n(d)$ and $\tilde{Z}$ as the same distribution as $\sup_{d\in \mathcal{D}_{B,s}} \xi_d$ and all assumptions are satisfied, based on the results of Section~\ref{sec:cvover}. For any $\varepsilon >0$ and $n$ large enough, we set $$\gamma_n= \frac{9 (bc\nu)^2}{\varepsilon} \frac{\ln^2 n}{a_{n,M_0,B} n^{1/3}}.$$
Thus, we have the following upper-bound
$$\begin{array}{ll}
\vspace{1em} \displaystyle \mathbb{P}\left(\left|\sup_{d\in \mathcal{D}_{B,s}} \mathcal{G}_n(d) - \sup_{d\in \mathcal{D}_{B,s}} \xi_d \right| > \sqrt{\varepsilon a_{n,M_0,B}} \right) &  \leq \mathbb{P} \left( |Z-\tilde Z| > \frac{3 bc\nu \ln n}{\gamma_n^{1/2} n^{1/6}} \right)\\
& \leq \mathbb{P}\left[
| Z - \tilde Z | >
\frac{b c\nu \ln n}{\gamma_n^{1/2}n^{1/4}} +
\frac{(b\sigma)^{1/2}(c\nu \ln n)^{3/4}}{\gamma_n^{1/2}n^{1/4}} + 
\frac{(b\sigma^2 (c\nu \ln n)^2)^{1/3}}{\gamma_n^{1/3}n^{1/6}}
\right].
\end{array}
$$
Using the Lemma~\ref{lem:empriproce}, we deduce that
$$\mathbb{P}\left(\left|\sup_{d\in \mathcal{D}_{B,s}} \mathcal{G}_n(d) - \sup_{d\in \mathcal{D}_{B,s}} \xi_d \right| > \sqrt{\varepsilon a_{n,M_0,B}} \right) \leq C\left(\gamma_n + \frac{\ln n}{n}\right),$$
and we conclude by noting that $\gamma_n$ goes to $0$ as $n$ goes to infinity.

\section{Additional numerical experiments} \label{sec:extraexpe}

\subsection{Investigating the impact of the rate of growing of $B$}
This section illustrates the impact of the rate of growing of $B$. We generate data from a mixture with three components and equal proportions ($\pi_k=1/3$). Each observations described by six variables. The density of $X_i$ given $Z_i$ is a product of univariate densities such that $X_{ij}=  \sum_{k=1}^K z_{ik}\delta_{kj} + \xi_{ij}$ where all the $\xi_{ij}$ are independent and where  $\delta_{11}=\delta_{12}=\delta_{23}=\delta_{24}=\delta_{35}=\delta_{36}=\tau$, while all remaining $\delta_{kj}=0$,  which implies that only the first six variables are relevant for clustering. Three distributions are considered for the $\xi_{ij}$ (standard Gaussian, Student with three degrees of freedom and Laplace) and the value of $\tau$ is defined to obtain a theoretical misclassification rate of $5\%$ ($\tau$ is equal to 1.94, 2.60 and 2.52 for the Gaussian, Student and Laplace distributions respectively). We consider different values of $B=[n^{1/s}]$ with $s=1,\ldots,10$. Note that the case $B=[n]$ does not satisfy the conditions stated by Assumption~\ref{ass:intervals}. Table~\ref{tab:bandwidth} indicates the empirical probability of selecting the true number of clusters obtained on 100 samples by considering that all the variables are relevant for clustering.
 Results illustrate that all the rates satisfying Assumption~\ref{ass:intervals} are asymptotically consistent. However, for small samples, results are improved by considering a small rate of growing ($[n^{1/6}]$ or slower). Note that the slow rates are not optimal for density estimation. However, the discretization is only used for model selection. Then, a kernel density estimation with an other bandwidth is performed for estimating the components of the selected model.

\begin{table}
\caption{Empirical probabilities of recovering the true number of clusters obtained with different numbers of bins by considering all the variables as relevant for clustering}\label{tab:bandwidth}
\fbox{%
\centering
\begin{tabular}{cccccccccccc}
  \hline
  distribution & sample size & $n^1$ & $n^{1/2}$ &$n^{1/3}$ &$n^{1/4}$ &$n^{1/5}$ &$n^{1/6}$ &$n^{1/7}$ &$n^{1/8}$ &$n^{1/9}$ & $ n^{1/10}$  \\ 
  \hline
  Gaussian & 100 & 0.00 & 0.00 & 0.04 & 0.74 & 0.74 & 0.83 & 0.83 & 0.83 & 0.83 & 0.83 \\ 
  & 250 & 0.00 & 0.00 & 0.94 & 1.00 & 1.00 & 1.00 & 1.00 & 1.00 & 1.00 & 1.00 \\ 
  & 500 & 0.00 & 0.00 & 1.00 & 1.00 & 1.00 & 1.00 & 1.00 & 1.00 & 1.00 & 1.00 \\ 
  & 1000 & 0.00 & 0.00 & 1.00 & 1.00 & 1.00 & 1.00 & 1.00 & 1.00 & 1.00 & 1.00 \\ 
  & 2000 & 0.00 & 0.00 & 1.00 & 1.00 & 1.00 & 1.00 & 1.00 & 1.00 & 1.00 & 1.00 \\ 
 & 5000 & 0.00 & 1.00 & 1.00 & 1.00 & 1.00 & 1.00 & 1.00 & 1.00 & 1.00 & 1.00 \\ 
  \hline
  Student & 100 & 0.00 & 0.00 & 0.01 & 0.85 & 0.85 & 0.90 & 0.90 & 0.90 & 0.90 & 0.90 \\ 
  & 250 & 0.00 & 0.00 & 0.96 & 1.00 & 1.00 & 1.00 & 1.00 & 1.00 & 1.00 & 1.00 \\ 
  & 500 & 0.00 & 0.00 & 1.00 & 1.00 & 1.00 & 1.00 & 1.00 & 1.00 & 1.00 & 1.00 \\ 
  & 1000 & 0.00 & 0.00 & 1.00 & 1.00 & 1.00 & 1.00 & 1.00 & 1.00 & 1.00 & 1.00 \\ 
 & 2000 & 0.00 & 0.00 & 1.00 & 1.00 & 1.00 & 1.00 & 1.00 & 1.00 & 1.00 & 1.00 \\ 
 & 5000 & 0.00 & 1.00 & 1.00 & 1.00 & 1.00 & 1.00 & 1.00 & 1.00 & 1.00 & 1.00 \\ 
  \hline
  Laplace & 100 & 0.00 & 0.00 & 0.05 & 0.81 & 0.81 & 0.89 & 0.89 & 0.89 & 0.89 & 0.89 \\ 
 & 250 & 0.00 & 0.00 & 0.99 & 1.00 & 1.00 & 1.00 & 1.00 & 1.00 & 1.00 & 1.00 \\ 
 & 500 & 0.00 & 0.00 & 1.00 & 1.00 & 1.00 & 1.00 & 1.00 & 1.00 & 1.00 & 1.00 \\ 
  & 1000 & 0.00 & 0.00 & 1.00 & 1.00 & 1.00 & 1.00 & 1.00 & 1.00 & 1.00 & 1.00 \\ 
  & 2000 & 0.00 & 0.02 & 1.00 & 1.00 & 1.00 & 1.00 & 1.00 & 1.00 & 1.00 & 1.00 \\ 
 & 5000 & 0.00 & 1.00 & 1.00 & 1.00 & 1.00 & 1.00 & 1.00 & 1.00 & 1.00 & 1.00 \\ 
     \hline
\end{tabular}
}
\end{table}

\subsection{Additional experiments for comparing the methods for a full model selection}
We consider the simulation setup presented in Section~\ref{sec:simucomparing} with two different theoretical misclassification rates.

\paragraph{Selection of the discriminative features}
In this experiment, we consider that the number of components is known. Thus, the task of model selection consists in determining the subset of the relevant variables. We consider the methods which automatically provide an estimator of the relevant variables (\emph{i.e.,} the proposed method, sparse $K$-means and VarSelLCM). Accuracy of this selection is measured by sensitivity (probability to detect as relevant a true discriminative variable) and specificity (probability to detect as irrelevant a true non discriminative variable). Table~\ref{tab:sen10} and \ref{tab:spe10} present the sensibility and the specificity obtained by the proposed approach and the parametric approach for a theoretical error rate of $10\%$.  Table~\ref{tab:sen15} and \ref{tab:spe15} present the sensibility and the specificity obtained by the proposed approach and the parametric approach for a theoretical error rate of $15\%$.
\begin{table}
\caption{Mean of the sensitivity (Sen.: $\text{card}(\widehat{\Omega}\cap\Omega)/6$) for the feature selection obtained by the proposed method (Proposed method), the parametric method (VarSelLCM) and the sparse K-means (Sparcl) on 100 replicates for each scenario with theoretical misclassification rate of $10\%$, when the number of components is known.} \label{tab:sen10}
\fbox{%
\centering
\begin{tabular}{ccccccccccc}
  \hline 
& & \multicolumn{3}{c}{Proposed method}& \multicolumn{3}{c}{VarSelLCM} & \multicolumn{3}{c}{Sparcl} \\
& & \multicolumn{3}{c}{$n$}& \multicolumn{3}{c}{$n$} & \multicolumn{3}{c}{$n$} \\
Component & $J$ &  100 & 250 & 500& 100 & 250 & 500& 100 & 250 & 500 \\
  \hline
  Gaussian & 20  & 0.42 & 0.92 & 1.00 & 0.70 & 1.00 & 1.00 & 0.88 & 0.94 & 0.97 \\ 
  & 50  & 0.24 & 0.43 & 1.00 & 0.36 & 0.99 & 1.00 & 0.83 & 0.98 & 1.00 \\ 
    & 100  & 0.11 & 0.14 & 0.86 & 0.15 & 0.85 & 1.00 & 0.63 & 0.96 & 0.99 \\ 
   Student & 20  & 0.59 & 0.98 & 1.00 & 0.18 & 0.26 & 0.36 & 0.76 & 0.81 & 0.88 \\ 
   & 50  & 0.33 & 0.66 & 1.00 & 0.11 & 0.17 & 0.23 & 0.60 & 0.81 & 0.84 \\ 
   & 100  & 0.23 & 0.17 & 0.99 & 0.11 & 0.18 & 0.20 & 0.42 & 0.70 & 0.82 \\ 
   Laplace & 20  & 0.65 & 0.98 & 1.00 & 0.24 & 0.54 & 0.86 & 0.85 & 0.86 & 0.88 \\ 
    & 50  & 0.42 & 0.73 & 1.00 & 0.09 & 0.18 & 0.33 & 0.78 & 0.94 & 0.95 \\ 
     & 100  & 0.22 & 0.24 & 0.98 & 0.03 & 0.03 & 0.08 & 0.63 & 0.94 & 0.98 \\ 
   \hline
\end{tabular}
}
\end{table}

\begin{table}
\caption{Mean of  the specificity (Spe.: $\text{card}(\widehat{\Omega}^c\cap\Omega^c)/(J-6)$) for the feature selection obtained by  the proposed method (Proposed method), the parametric method (VarSelLCM) and the sparse K-means (Sparcl)  on 100 replicates for each scenario with theoretical misclassification rate of $10\%$, when the number of components is known.} \label{tab:spe10}
\fbox{%
\centering
\begin{tabular}{ccccccccccc}
  \hline 
& & \multicolumn{3}{c}{Proposed method}& \multicolumn{3}{c}{VarSelLCM} & \multicolumn{3}{c}{Sparcl} \\
& & \multicolumn{3}{c}{$n$}& \multicolumn{3}{c}{$n$} & \multicolumn{3}{c}{$n$} \\
Component & $J$ &  100 & 250 & 500& 100 & 250 & 500& 100 & 250 & 500 \\
  \hline
  Gaussian & 20  & 0.98 & 1.00 & 1.00 & 0.99 & 1.00 & 1.00 & 0.78 & 0.69 & 0.57 \\ 
  & 50 & 0.97 & 1.00 & 1.00 & 1.00 & 1.00 & 1.00 & 0.67 & 0.43 & 0.32 \\ 
  & 100 & 0.98 & 1.00 & 1.00 & 1.00 & 1.00 & 1.00 & 0.69 & 0.42 & 0.22 \\ 
 Student & 20 & 0.97 & 1.00 & 1.00 & 0.68 & 0.57 & 0.44 & 0.67 & 0.56 & 0.54 \\ 
  & 50 & 0.97 & 1.00 & 1.00 & 0.74 & 0.67 & 0.59 & 0.70 & 0.73 & 0.58 \\ 
  & 100  & 0.98 & 1.00 & 1.00 & 0.76 & 0.70 & 0.62 & 0.82 & 0.75 & 0.68 \\ 
 Laplace & 20  & 0.98 & 1.00 & 1.00 & 0.85 & 0.85 & 0.89 & 0.74 & 0.83 & 0.82 \\ 
   & 50  & 0.98 & 1.00 & 1.00 & 0.90 & 0.90 & 0.87 & 0.63 & 0.77 & 0.73 \\ 
   & 100  & 0.99 & 1.00 & 1.00 & 0.91 & 0.93 & 0.91 & 0.71 & 0.49 & 0.60 \\
   \hline
\end{tabular}
}
\end{table}

\begin{table}
\caption{Mean of the sensitivity (Sen.: $\text{card}(\widehat{\Omega}\cap\Omega)/6$) for the feature selection obtained by the proposed method (Proposed method), the parametric method (VarSelLCM) and the sparse K-means (Sparcl) on 100 replicates for each scenario with theoretical misclassification rate of $15\%$, when the number of components is known.}\label{tab:sen15}
\fbox{%
\centering
\begin{tabular}{ccccccccccc}
  \hline 
& & \multicolumn{3}{c}{Proposed method}& \multicolumn{3}{c}{VarSelLCM} & \multicolumn{3}{c}{Sparcl} \\
& & \multicolumn{3}{c}{$n$}& \multicolumn{3}{c}{$n$} & \multicolumn{3}{c}{$n$} \\
Component & $J$ &  100 & 250 & 500& 100 & 250 & 500& 100 & 250 & 500 \\
  \hline
 Gaussian & 20  & 0.26 & 0.43 & 0.96 & 0.24 & 0.87 & 1.00 & 0.88 & 0.97 & 0.99 \\ 
   & 50 & 0.12 & 0.05 & 0.72 & 0.07 & 0.63 & 1.00 & 0.70 & 0.99 & 1.00 \\ 
   & 100  & 0.04 & 0.01 & 0.36 & 0.03 & 0.31 & 0.95 & 0.47 & 0.94 & 1.00 \\ 
  Student & 20  & 0.34 & 0.60 & 1.00 & 0.18 & 0.30 & 0.39 & 0.70 & 0.89 & 0.92 \\ 
   & 50  & 0.10 & 0.17 & 0.89 & 0.13 & 0.20 & 0.29 & 0.51 & 0.80 & 0.92 \\ 
    & 100 0 & 0.09 & 0.01 & 0.50 & 0.13 & 0.22 & 0.23 & 0.31 & 0.57 & 0.80 \\ 
 Laplace & 20  & 0.31 & 0.63 & 1.00 & 0.10 & 0.09 & 0.11 & 0.82 & 0.90 & 0.90 \\ 
  & 50  & 0.20 & 0.18 & 0.91 & 0.05 & 0.02 & 0.03 & 0.64 & 0.96 & 0.97 \\ 
  & 100  & 0.10 & 0.03 & 0.55 & 0.04 & 0.01 & 0.02 & 0.46 & 0.87 & 0.99 \\ 
   \hline
\end{tabular}
}
\end{table}

\begin{table}
\caption{Mean of  the specificity (Spe.: $\text{card}(\widehat{\Omega}^c\cap\Omega^c)/(J-6)$) for the feature selection obtained by  the proposed method (Proposed method), the paramteric method (VarSelLCM) and the sparse K-means (Sparcl)  on 100 replicates for each scenario with theoretical misclassification rate of $15\%$, when the number of components is known.} \label{tab:spe15}
\fbox{%
\centering
\begin{tabular}{ccccccccccc}
  \hline 
& & \multicolumn{3}{c}{Proposed method}& \multicolumn{3}{c}{VarSelLCM} & \multicolumn{3}{c}{Sparcl} \\
& & \multicolumn{3}{c}{$n$}& \multicolumn{3}{c}{$n$} & \multicolumn{3}{c}{$n$} \\
Component & $J$ &  100 & 250 & 500& 100 & 250 & 500& 100 & 250 & 500 \\
  \hline
    Gaussian & 20 & 0.97 & 0.99 & 1.00 & 0.97 & 1.00 & 1.00 & 0.62 & 0.56 & 0.40 \\ 
   & 50 & 0.97 & 1.00 & 1.00 & 0.99 & 1.00 & 1.00 & 0.56 & 0.33 & 0.26 \\ 
   & 100  & 0.99 & 1.00 & 1.00 & 1.00 & 1.00 & 1.00 & 0.69 & 0.27 & 0.15 \\ 
 Student & 20  & 0.96 & 1.00 & 1.00 & 0.69 & 0.57 & 0.44 & 0.55 & 0.43 & 0.39 \\ 
  & 50 & 0.98 & 1.00 & 1.00 & 0.74 & 0.67 & 0.58 & 0.66 & 0.52 & 0.40 \\ 
   & 100  & 0.98 & 1.00 & 1.00 & 0.76 & 0.69 & 0.63 & 0.81 & 0.63 & 0.52 \\ 
  Laplace & 20  & 0.97 & 0.99 & 1.00 & 0.82 & 0.75 & 0.65 & 0.55 & 0.64 & 0.58 \\ 
   & 50 & 0.97 & 1.00 & 1.00 & 0.89 & 0.88 & 0.83 & 0.60 & 0.45 & 0.54 \\ 
   & 100  & 0.98 & 1.00 & 1.00 & 0.91 & 0.93 & 0.91 & 0.70 & 0.34 & 0.30 \\ 
   \hline
\end{tabular}
}
\end{table}

\paragraph{Full model selection}
We now compare both non-parametric and parametric approaches on their performances for full model selection. Table~\ref{tab:selectK10} and  Table~\ref{tab:selectK15} present the statistics on the number of components selected by both approaches, when the theoretical misclassification rate is $10\%$ and $15\%$ respectively.

\begin{table}
\caption{Probability to select the true number of components (Tr.) and to overestimate it (Ov.) obtained by the proposed method (Proposed method) and the parametric method (VarSelLCM) on 100 replicates for each scenario  with theoretical misclassification rate of $10\%$, by performing a selection of  the variables.} \label{tab:selectK10}
\fbox{%
\centering
\begin{tabular}{ccccccccccccccc}
  \hline  
Component & $J$ & \multicolumn{6}{c}{Proposed method}& \multicolumn{6}{c}{VarSelLCM} \\
& & \multicolumn{2}{c}{$n=100$}& \multicolumn{2}{c}{$n=250$}& \multicolumn{2}{c}{$n=500$}& \multicolumn{2}{c}{$n=100$}& \multicolumn{2}{c}{$n=250$}& \multicolumn{2}{c}{$n=500$}\\
& & Tr. & Ov. & Tr. & Ov.& Tr. & Ov.& Tr. & Ov.& Tr. & Ov.& Tr. & Ov. \\ 
\hline
 Gaussian & 20  & 0.07 & 0.00 & 0.69 & 0.00 & 1.00 & 0.00 & 0.21 & 0.00 & 1.00 & 0.00 & 1.00 & 0.00 \\ 
    & 50  & 0.00 & 0.00 & 0.15 & 0.00 & 0.98 & 0.01 & 0.06 & 0.00 & 1.00 & 0.00 & 1.00 & 0.00 \\ 
     & 100  & 0.01 & 0.00 & 0.02 & 0.00 & 0.73 & 0.00 & 0.01 & 0.00 & 0.81 & 0.01 & 1.00 & 0.00 \\ 
    Student & 20  & 0.13 & 0.00 & 0.89 & 0.00 & 1.00 & 0.00 & 0.62 & 0.08 & 0.55 & 0.42 & 0.04 & 0.96 \\ 
     & 50  & 0.06 & 0.00 & 0.45 & 0.00 & 1.00 & 0.00 & 0.80 & 0.16 & 0.41 & 0.59 & 0.03 & 0.97 \\ 
     & 100  & 0.00 & 0.00 & 0.01 & 0.00 & 0.95 & 0.00 & 0.68 & 0.32 & 0.13 & 0.87 & 0.00 & 1.00 \\ 
    Laplace & 20  & 0.17 & 0.00 & 0.96 & 0.00 & 1.00 & 0.00 & 0.18 & 0.01 & 0.58 & 0.24 & 0.05 & 0.95 \\ 
    & 50  & 0.07 & 0.00 & 0.58 & 0.00 & 1.00 & 0.00 & 0.10 & 0.00 & 0.22 & 0.04 & 0.35 & 0.59 \\ 
   & 100  & 0.02 & 0.00 & 0.08 & 0.00 & 0.99 & 0.00 & 0.03 & 0.00 & 0.03 & 0.00 & 0.28 & 0.20 \\ 
   \hline
\end{tabular}
}
\end{table}

\begin{table}
\caption{Probability to select the true number of components (Tr.) and to overestimate it (Ov.) obtained by the proposed method (Proposed method) and the parametric method (VarSelLCM) on 100 replicates for each scenario  with theoretical misclassification rate of $15\%$, by performing a selection of  the variables.\label{tab:selectK15}
}
\fbox{%
\centering
\begin{tabular}{ccccccccccccccc}
  \hline  
Component & $J$ & \multicolumn{6}{c}{Proposed method}& \multicolumn{6}{c}{VarSelLCM} \\
& & \multicolumn{2}{c}{$n=100$}& \multicolumn{2}{c}{$n=250$}& \multicolumn{2}{c}{$n=500$}& \multicolumn{2}{c}{$n=100$}& \multicolumn{2}{c}{$n=250$}& \multicolumn{2}{c}{$n=500$}\\
& & Tr. & Ov. & Tr. & Ov.& Tr. & Ov.& Tr. & Ov.& Tr. & Ov.& Tr. & Ov. \\ 
  \hline
 Gaussian & 20  & 0.00 & 0.00 & 0.02 & 0.00 & 0.86 & 0.00 & 0.04 & 0.00 & 0.50 & 0.00 & 1.00 & 0.00 \\ 
   & 50  & 0.00 & 0.00 & 0.00 & 0.00 & 0.48 & 0.00 & 0.02 & 0.00 & 0.30 & 0.00 & 0.99 & 0.01 \\ 
  & 100  & 0.00 & 0.00 & 0.00 & 0.00 & 0.13 & 0.00 & 0.00 & 0.00 & 0.08 & 0.00 & 0.82 & 0.02 \\ 
 Student & 20 & 0.00 & 0.00 & 0.02 & 0.00 & 0.86 & 0.00 & 0.04 & 0.00 & 0.50 & 0.00 & 1.00 & 0.00 \\ 
   & 50  & 0.00 & 0.00 & 0.00 & 0.00 & 0.48 & 0.00 & 0.02 & 0.00 & 0.30 & 0.00 & 0.99 & 0.01 \\ 
   & 100  & 0.00 & 0.00 & 0.00 & 0.00 & 0.13 & 0.00 & 0.00 & 0.00 & 0.08 & 0.00 & 0.82 & 0.02 \\ 
 Laplace & 20  & 0.00 & 0.00 & 0.02 & 0.00 & 0.86 & 0.00 & 0.04 & 0.00 & 0.50 & 0.00 & 1.00 & 0.00 \\ 
  & 50 & 0.00 & 0.00 & 0.00 & 0.00 & 0.48 & 0.00 & 0.02 & 0.00 & 0.30 & 0.00 & 0.99 & 0.01 \\ 
   & 100  & 0.00 & 0.00 & 0.00 & 0.00 & 0.13 & 0.00 & 0.00 & 0.00 & 0.08 & 0.00 & 0.82 & 0.02 \\ 
   \hline
\end{tabular}
}
\end{table}

Table~\ref{tab:simu2feature10} and Table~\ref{tab:simu2feature15} present the sensitivity and the specificity for feature selection obtained by both approaches when the number of components is also estimated, under a theoretical misclassification rate of $10\%$ and $15\%$ respectively.
\begin{table}
\caption{Mean of the sensitivity (Sen.: $\text{card}(\widehat{\Omega}\cap\Omega)/6$) and the specificity (Spe.: $\text{card}(\widehat{\Omega}^c\cap\Omega^c)/(J-6)$) for the feature selection obtained by the proposed method (Proposed method) and the parametric method (VarSelLCM) on 100 replicates for each scenario with theoretical misclassification rate of $10\%$, when the number of components also is estimated.}
\label{tab:simu2feature10}
\fbox{%
\centering
\begin{tabular}{ccccccccccccccc}
  \hline  
Component & $J$ & \multicolumn{6}{c}{Proposed method}& \multicolumn{6}{c}{VarSelLCM} \\
& & \multicolumn{2}{c}{$n=100$}& \multicolumn{2}{c}{$n=250$}& \multicolumn{2}{c}{$n=500$}& \multicolumn{2}{c}{$n=100$}& \multicolumn{2}{c}{$n=250$}& \multicolumn{2}{c}{$n=500$}\\
& &  Sen. & Spe. & Sen. & Spe. & Sen. & Spe. & Sen. & Spe. & Sen. & Spe. & Sen. & Sep. \\ 
  \hline
 Gaussian & 20  & 0.88 & 0.61 & 1.00 & 0.90 & 1.00 & 1.00 & 0.97 & 0.69 & 1.00 & 1.00 & 1.00 & 1.00 \\ 
    & 50  & 0.93 & 0.46 & 0.99 & 0.72 & 1.00 & 1.00 & 0.98 & 0.59 & 1.00 & 1.00 & 1.00 & 1.00 \\ 
   & 100  & 0.92 & 0.31 & 0.94 & 0.67 & 1.00 & 0.94 & 0.95 & 0.42 & 1.00 & 0.96 & 1.00 & 1.00 \\ 
  Student & 20  & 0.91 & 0.67 & 1.00 & 0.97 & 1.00 & 1.00 & 0.66 & 0.18 & 0.56 & 0.28 & 0.41 & 0.42 \\ 
     & 50  & 0.93 & 0.51 & 1.00 & 0.82 & 1.00 & 1.00 & 0.74 & 0.10 & 0.68 & 0.17 & 0.58 & 0.24 \\ 
     & 100  & 0.92 & 0.39 & 1.00 & 0.71 & 1.00 & 0.99 & 0.77 & 0.12 & 0.72 & 0.16 & 0.67 & 0.16 \\ 
    Laplace & 20 & 0.94 & 0.70 & 1.00 & 0.98 & 1.00 & 1.00 & 0.79 & 0.23 & 0.85 & 0.56 & 0.88 & 0.93 \\ 
     & 50 & 0.94 & 0.55 & 1.00 & 0.88 & 1.00 & 1.00 & 0.82 & 0.11 & 0.86 & 0.16 & 0.89 & 0.50 \\ 
    & 100  & 0.93 & 0.44 & 1.00 & 0.72 & 1.00 & 1.00 & 0.84 & 0.07 & 0.88 & 0.06 & 0.90 & 0.08 \\ 
   \hline
\end{tabular}
}
\end{table}

\begin{table}
\caption{Mean of the sensitivity (Sen.: $\text{card}(\widehat{\Omega}\cap\Omega)/6$) and the specificity (Spe.: $\text{card}(\widehat{\Omega}^c\cap\Omega^c)/(J-6)$) for the feature selection obtained by the proposed method (Proposed method) and the parametric method (VarSelLCM) on 100 replicates for each scenario with theoretical misclassification rate of $15\%$, when the number of components also is estimated.}
\label{tab:simu2feature15}
\fbox{%
\centering
\begin{tabular}{ccccccccccccccc}
  \hline  
Component & $J$ & \multicolumn{6}{c}{Proposed method}& \multicolumn{6}{c}{VarSelLCM} \\
& & \multicolumn{2}{c}{$n=100$}& \multicolumn{2}{c}{$n=250$}& \multicolumn{2}{c}{$n=500$}& \multicolumn{2}{c}{$n=100$}& \multicolumn{2}{c}{$n=250$}& \multicolumn{2}{c}{$n=500$}\\
& &  Sen. & Spe. & Sen. & Spe. & Sen. & Spe. & Sen. & Spe. & Sen. & Spe. & Sen. & Sep. \\ 
  \hline
  Gaussian & 20 & 0.67 & 0.66 & 0.80 & 0.68 & 1.00 & 0.95 & 0.90 & 0.48 & 1.00 & 0.84 & 1.00 & 1.00 \\ 
    & 50  & 0.90 & 0.29 & 0.67 & 0.69 & 1.00 & 0.82 & 0.93 & 0.33 & 0.99 & 0.76 & 1.00 & 1.00 \\ 
    & 100  & 0.91 & 0.16 & 0.56 & 0.70 & 0.99 & 0.70 & 0.95 & 0.25 & 0.99 & 0.63 & 1.00 & 0.96 \\ 
   Student & 20  & 0.67 & 0.66 & 0.80 & 0.68 & 1.00 & 0.95 & 0.90 & 0.48 & 1.00 & 0.84 & 1.00 & 1.00 \\ 
    & 50  & 0.90 & 0.29 & 0.67 & 0.69 & 1.00 & 0.82 & 0.93 & 0.33 & 0.99 & 0.76 & 1.00 & 1.00 \\ 
   & 100  & 0.91 & 0.16 & 0.56 & 0.70 & 0.99 & 0.70 & 0.95 & 0.25 & 0.99 & 0.63 & 1.00 & 0.96 \\ 
   Laplace & 20  & 0.67 & 0.66 & 0.80 & 0.68 & 1.00 & 0.95 & 0.90 & 0.48 & 1.00 & 0.84 & 1.00 & 1.00 \\ 
    & 50  & 0.90 & 0.29 & 0.67 & 0.69 & 1.00 & 0.82 & 0.93 & 0.33 & 0.99 & 0.76 & 1.00 & 1.00 \\ 
    & 100  & 0.91 & 0.16 & 0.56 & 0.70 & 0.99 & 0.70 & 0.95 & 0.25 & 0.99 & 0.63 & 1.00 & 0.96 \\
   \hline
\end{tabular}
}
\end{table}

\begin{figure}\label{fig:resultsfeatureselection10}
\centering \includegraphics[scale=0.36]{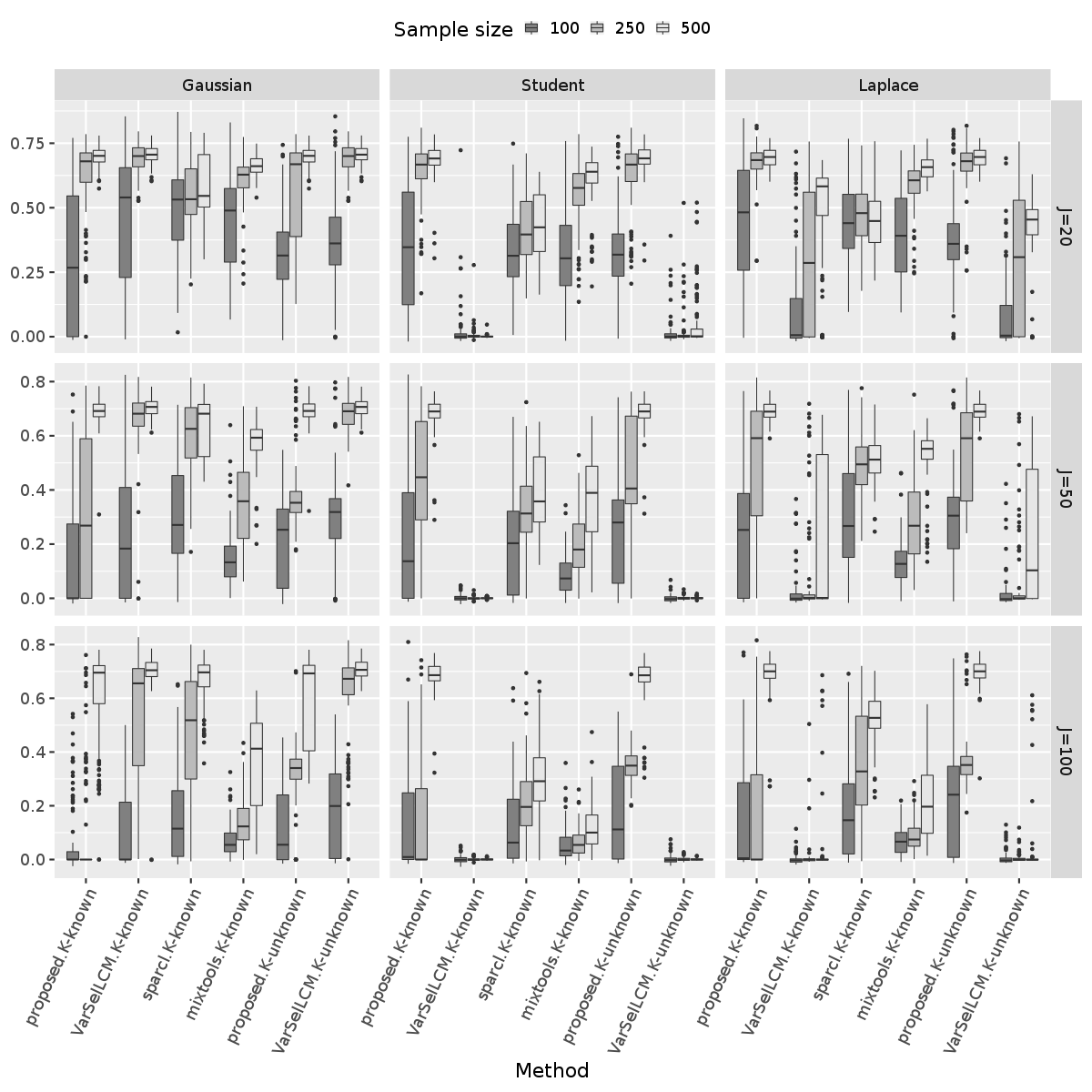}
\caption{Boxplot of the Adjusted Rand Index (ARI) obtained on the resulting partition when feature selection is performed with the true number of components by the proposed method (proposed.K-known) and by the parametric method (VarSelLCM.K-known), by the sparse K-means (Sparcl.K-known) and by the model considering all the variables as relevant components (mixtools.K-known) and when the full model selection (feature selection and estimation of the number of components) is achieved by the proposed approach (proposed.K-unknown) and the parametric approach (VarSelLCM.K-unknown). Data are generated with theoretical misclassification rate of $10\%$.}
\end{figure}

\begin{figure}\label{fig:resultsfeatureselection15}
\centering \includegraphics[scale=0.36]{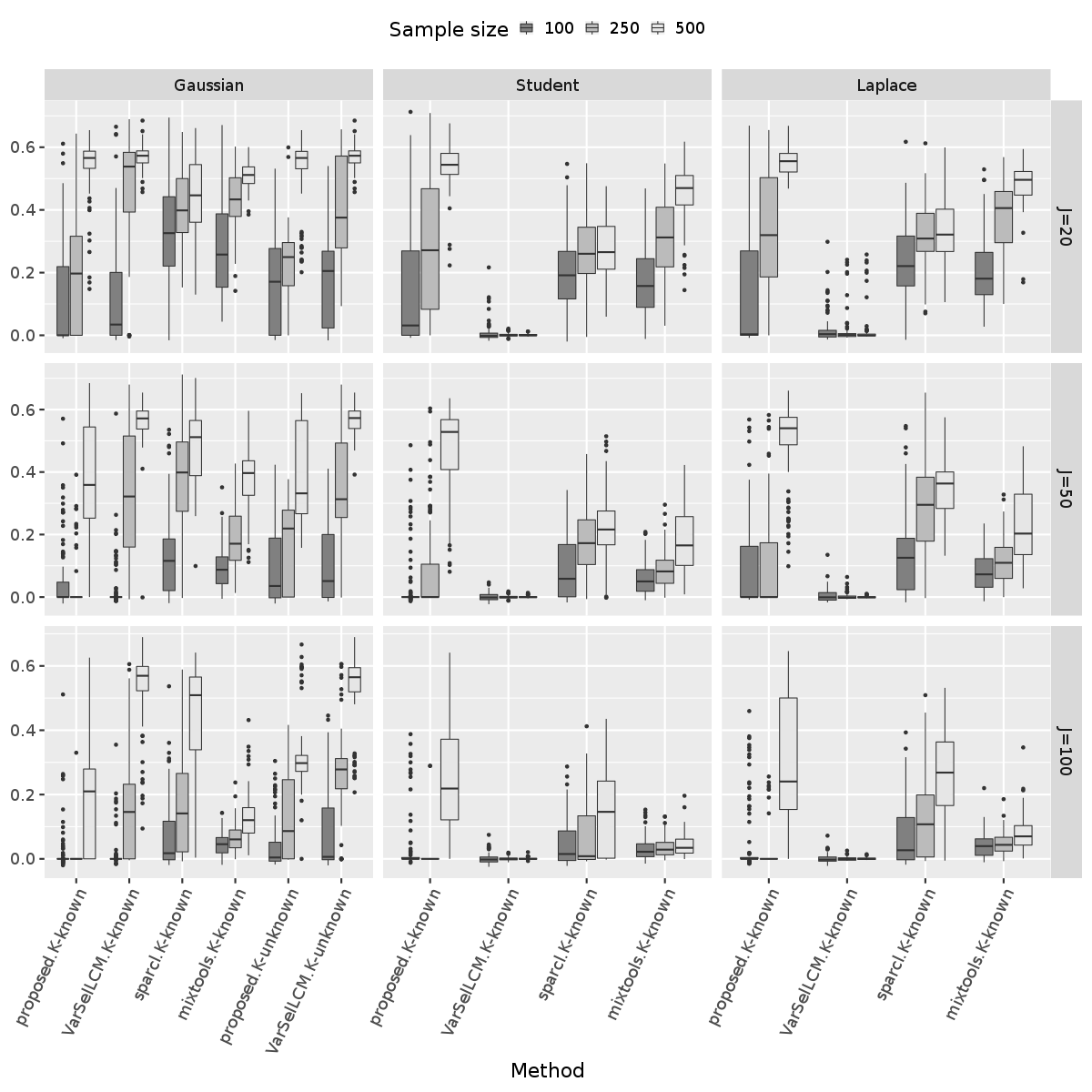}
\caption{Boxplot of the Adjusted Rand Index (ARI) obtained on the resulting partition when feature selection is performed with the true number of components by the proposed method (proposed.K-known) and by the parametric method (VarSelLCM.K-known), by the sparse K-means (Sparcl.K-known) and by the model considering all the variables as relevant components (mixtools.K-known) and when the full model selection (feature selection and estimation of the number of components) is achieved by the proposed approach (proposed.K-unknown) and the parametric approach (VarSelLCM.K-unknown). Data are generated with theoretical misclassification rate of $15\%$.}
\end{figure}

We are now interested in investigating the accuracy of the estimated partition. 
Figure~\ref{fig:resultsfeatureselection10} and Figure~\ref{fig:resultsfeatureselection15} present the distribution of the ARI obtained when $K_0$ is known or estimated,  for a theoretical misclassification rate of $10\%$ and $15\%$ respectively.

 \bibliographystyle{apalike}
\bibliography{biblio}